\documentclass[10pt]{amsart}
\usepackage{amsmath, amscd, amsthm} 
\usepackage{amssymb, amsfonts} 
\usepackage{spseq} 
\usepackage[all]{xy} 
\begin{document}
\title[Equivariant Invariants of Real Varieties]{Equivariant Semi-topological Invariants, Atiyah's $KR$-theory, and Real Algebraic Cycles}
\author{Jeremiah Heller}
\email{heller@math.northwestern.edu}
\author{Mircea Voineagu}
\email{voineagu@usc.edu, mircea.voineagu@ipmu.jp}

\begin{abstract}
We establish an Atiyah-Hirzebruch type spectral sequence relating real morphic cohomology and real semi-topological $K$-theory and prove it to be compatible with the Atiyah-Hirzebruch spectral sequence relating Bredon cohomology and Atiyah's KR-theory constructed by Dugger. An equivariant  and a real version of Suslin's conjecture on morphic cohomology are formulated, proved to come from the complex version of Suslin conjecture and verified for certain real varieties. In conjunction with the spectral sequences constructed here this allows the computation of the real semi-topological $K$-theory of some real varieties. As another application of this spectral sequence we give an alternate proof of the Lichtenbaum-Quillen conjecture over $\R$, extending an earlier proof of Karoubi and Weibel.
\end{abstract}

\maketitle
\tableofcontents
\section{Introduction}

A semi-topological Atiyah-Hirzebruch spectral sequence relating the morphic cohomology and the semi-topological $K$-theory of  smooth complex variety $X$ is constructed is constructed in \cite{FHW:sst}. There is a comparison map from the motivic Atiyah-Hirzebruch spectral sequence to the semi-topological spectral sequence and a map from the semi-topological spectral sequence to the topological Atiyah-Hirzebruch spectral sequence for connective complex $K$-theory. A first goal of this paper is to fill in this picture for real varieties. A second goal is to compute equivariant morphic cohomology and real semi-topological $K$-theory for some real varieties.

In Section \ref{semiss} we establish the semi-topological Atiyah-Hirzebruch spectral sequence and show that it receives a map from the motivic spectral sequence and maps to the spectral sequence constructed by Dugger in \cite[Corollary 1.3]{dug:kr}, 
\begin{equation*}
 \begin{CD}
  E_{2}^{p,q}(alg)= \H^{p-q}_{\mcal{M}}(X ; A(-q)) \Longrightarrow K_{-p-q}(X;A) \\
@VVV \\
E_{2}^{p,q}(sst) = L^{-q}H\R^{p,-q}(X; A) \Longrightarrow K\R^{sst}_{-p-q}(X;A) \\
@VVV \\
E_{2}^{p,q}(top) = H^{p,-q}(X(\C) ; \underline{A}) \Longrightarrow kr^{p+q}(X(\C);A).
 \end{CD}
\end{equation*}

To construct the real semi-topological spectral sequence we apply Friedlander-Haesemeyer-Walker's method \cite{FHW:sst}, with straightforward modifications for real varieties. As their method makes clear, in order to produce the semi-topological spectral sequence one may use as input any model for the motivic spectral sequence. However, constructing the companion topological spectral sequence and identifying the equivariant homotopy type of the spectra appearing in its construction requires more care. Presumably this can be carried out with any of the available models, however we find it convenient to make use of Grayson's model (as described in \cite{Walk:thes, Walker:adams}). Suslin \cite{Suslin:GraySS} has proved that the motivic cohomology appears as the $E_{2}$-term of Grayson's spectral sequence. Moreover it is widely believed that Friedlander-Suslin's spectral sequence \cite{FS:AHSS}, Levine's spectral sequence \cite{Levine:ss}, and Grayson's spectral sequence are in fact all the same spectral sequence. The advantage that this model has for us is that the spectra in this tower arise from homotopy group completions of simplicial sets of algebraic maps to certain quot schemes. This geometric description makes  it straightforward to prove Proposition \ref{gres}, which is one of the main technical tools we use to identify the form of our topological spectral sequence. The topological spectral sequence we construct here has the same form as the spectral sequence constructed by Dugger in \cite[Corollary 1.3]{dug:kr} and indeed in Proposition \ref{dugagr} we prove that these two spectral sequences are the same. 

Suslin's conjecture on morphic cohomology of a smooth complex variety predicts that the cycle map $L^{q}H^{n}(X) \to H_{sing}^{n}(X(\C);\Z)$ is an isomorphism for $n\leq q$ and an injection for $n = q+1$. We formulate appropriate analogues for the real morphic cohomology and the equivariant morphic cohomology of a smooth real variety (see Section \ref{conj} for precise statements). The real morphic cohomology groups of real varieties and the morphic cohomology of complex varieties are particular instances of equivariant morphic cohomology, and so at first glance the equivariant conjecture would appear to be stronger than the rest of these conjectures. However, as we show in Theorem \ref{equiv} all three of these conjectures are equivalent. As a consequence we verify the equivariant Suslin conjecture in codimension one for all smooth real varieties and give a complete computation of real morphic cohomology in codimension one (see Theorem \ref{div}), generalizing a result of Teh \cite{Teh:HT}. For certain real varieties, such as geometrically rationally connected threefolds and generic cubics of small dimension we verify the equivariant Suslin conjecture in all codimensions.

Moreover in Section \ref{conj}, we discuss the real algebraic equivalence relation on real cycles (as defined by \cite{FW:real}). We prove that it satisfies a Bloch-Ogus formula for Bredon cohomology (and also Borel cohomology) as well as a number of other  good properties satisfied by the usual algebraic equivalence over complex algebraic cycles (see Theorem \ref{cod2} and Theorem \ref{div}). Together with other techniques, these allow us to prove a stronger result than is predicted in general by the real Suslin's conjecture in some cases (e.g. geometrically rationally connected threefolds), extending results from the complex case of Voineagu \cite{Voin:1}. In particular we obtain that all of the equivariant morphic cohomology groups are finitely generated in these cases.

Using the spectral sequences constructed in Section \ref{semiss} in conjunction with our computations of the real morphic cohomology, we compute the real semi-topological $K$-theory of some real varieties.  As an additional application of the spectral sequences constructed here, in conjunction with Voevodsky's verification of the Milnor conjecture, we give an alternate proof of the $2$-adic Lichtenbaum-Quillen over $\R$. This was proved \cite{KW:real} (in a slightly weaker form) and in general (for any real closed field) in \cite{RO:LQ}. One would like to use Voevodsky's solution of Milnor conjecture \cite{Voev:miln} together with a comparison of Atiyah-Hirzebruch spectral sequences in order to prove this conjecture (indeed this method works just fine over non-real closed fields). The difficulty which has to be worked around in \cite{KW:real} and \cite{RO:LQ} is that if $X$ has a real point then its etale $2$-cohomological dimension is infinite and the etale Atiyah-Hirzebruch spectral sequence will not converge. The proof we give avoids this problem by using the Atiyah-Hirzebruch spectral sequence based on the Bredon cohomology of the $\Z/2$-space $X(\C)$ (which is convergent). We view this application as an example of the fact that Bredon cohomology is a more natural and useful target for Chow groups of real algebraic cycles than Borel cohomology.

Finally in Section \ref{comppd} we verify the expected compatibility between the duality proved in \cite{HV:VT} relating equivariant morphic cohomology and dos Santos' equivariant Lawson homology groups of a smooth projective variety $X$ and the Poincare duality relating the Bredon cohomology and homology of the space $X(\C)$. An important consequence of this compatibility is Corollary \ref{geqd} which asserts that in weights larger than $\dim U$, the cycle map from the equivariant morphic cohomology of $U$ to the Bredon cohomology of $U$ is an isomorphism for any smooth quasi-projective real variety $U$. This identification is an important ingredient in the previously mentioned applications of the spectral sequences. Another application of this compatibility was given in \cite[Corollary 5.14]{HV:VT}.

The  techniques and computations presented in this paper are a sequel to our previous work on real algebraic cycles \cite{HV:VT} and to the paper \cite{FHW:sst} of Friedlander-Haesemeyer-Walker which had an important  influence on this work.

{\bf{Acknowledgements}} We were told by E. Friedlander that a version of the spectral sequence constructed in section \ref{semiss} of this paper was also constructed independently by E. Friedlander, M. Walker and C. Haesemeyer. We thank him for informing us about their work, for his interest in our paper, and his encouragements. We also thank Pedro dos Santos and Jens Hornbostel for interest in this paper and Mark Walker for useful discussions.

\subsection{Notation}
By a $k$-variety we mean a reduced, separated quasi-projective scheme of finite type over a field $k$. We write $Sch/k$ for the category of $k$-varieties  and  $Sm/k$ for the subcategory of smooth quasi-projective $k$-varieties.  Unless otherwise specified, $G$ denotes the group $Gal(\C/\R)=\Z/2$.

If $X$ is a real variety then the set $\Hom{\R}{\C}{X}$ is denoted as $X(\C)$. This set agrees with the set of closed points $X_{\C}(\C)$ of the complexification and is given an analytic topology via this identification. Note that when $X$ is a complex variety viewed as a real variety this notation is potentially ambiguous but it will always be clear in such situations which is meant. Via the action of complex conjugation the space $X(\C)$ is a $G$-space. 

By a spectrum we mean a spectrum in the sense of Bousfield-Friedlander, which consists of a sequence $A = (A_{1}, A_{2}, \ldots)$ of simplicial sets together with connecting maps $\sigma_{i}:S^{1}\wedge A_{i} \to A_{i+1}$. A naive $G$-spectrum is a $G$-object in the category of spectra.  

A map of $G$-spaces (or naive $G$-spectra) $f:A\to B$ is a $G$-weak equivalence provided both $A\to B$ and $A^{G}\to B^{G}$ are weak equivalences.

Let $V$ be a real representation of $G$, write $S^{V}$ for the one-point compactification of $V$. The $V$th homotopy group of a based $G$-space $X$ is
\begin{equation*}
 \pi_{V} X = [S^{V}, X]_{G} ,
\end{equation*}
(where $[-,-]_{G}$ denotes maps in the based $G$-homotopy category).
When  $V = \R^{p,q}$ we use the notation
$\pi_{p,q}X = \pi_{\R^{p,q}}X$ where $\R^{p,q}$ is $\R^{p+q}$ with trivial action on the first $p$-components and is multiplication by $-1$ on the last $q$-components. In particular $\pi_{k,0}X = \pi_{k}X^{G}$. We use a similar convention for the $RO(G)$-graded Bredon cohomology with coefficients in a Mackey functor $\underline{M}$ and write  $H^{p,q}(X;\underline{M}) = H^{\R^{p,q}}(X;\underline{M})$. For more details on equivariant homotopy and cohomology see \cite{May:equi}.

\section{Atiyah-Hirzebruch Spectral Sequences}\label{semiss}
In this section we construct the semi-topological Atiyah-Hirzebruch spectral sequence for real varieties following the construction in \cite{FHW:sst} for complex varieties. We show that this spectral sequence receives a map from the motivic Atiyah-Hirzebruch spectral sequence and maps to the Atiyah-Hirzebruch spectral sequence for Bredon cohomology and real $K$-theory constructed by Dugger. 

\subsection{Motivic spectral sequence}\label{motss}
There are several available constructions of a motivic Atiyah-Hirzebruch spectral sequence. For our purposes we find it convenient to make use of Grayson's construction \cite{Gray:wt} (as reinterpreted by Walker \cite{Walk:thes,Walker:adams}).

Let $f$ be a polynomial (integer valued with rational coefficients), $\mcal{E}$ a coherent sheaf on a quasi-projective variety $Y$, and $Y\to S$ a morphism of quasi-projective varieties. The quot scheme $\Quot^{f}_{\mcal{E}/Y/S}$ represents the functor which sends $X$ to the 
collection of quotient objects $[\pi_{X}^*\mcal{E}\epi \mcal{F}]$ (where $\pi_{X}:X\times_{S}Y\to Y$ is the projection) which are flat over $X$ and have support proper over $X$ with Hilbert polynomial $f$.

\begin{definition}\label{Gn}
For an integer $n$ define $G_{Y}(n) \subseteq \coprod_{\deg(f) = 0} \Quot^{f}_{\mcal{O}^{n}_{Y}/Y/\spec k}$ to be the subfunctor 
whose value on $X$ is the collection of quotient objects $[\mcal{O}^{n}_{X\times Y}\twoheadrightarrow \mcal{F}]$ which satisfy the conditions
\begin{enumerate}
 \item $\mcal{F}$ is flat over $X$ and $Supp(\mcal{F})$ is finite over $X$, 
 \item The induced map $\mcal{O}^{n}_{X} \to \pi_{*}\mcal{O}^{n}_{X\times Y} \to \pi_{*}\mcal{F}$ is surjective.
\end{enumerate}
\end{definition}
By \cite[Lemma 2.2]{Walk:Thom} this functor is represented by a variety which we also denote $G_{Y}(n)$. In fact $G_{Y}(n)\subseteq \coprod^{n}_{\deg(f) = 0} \Quot^{f}_{\mcal{O}^{n}_{Y}/Y/\spec k} $ is an open subvariety. In \cite{GW:geom} Grayson-Walker use  $K_{Y}(n)=  \coprod_{\deg(f) = 0} \Quot^{f}_{\mcal{O}^{n}_{Y}/Y/\spec k}$ rather than $G_{Y}(n)$. Using $G_{Y}(n)$ rather than $K_{Y}(n)$ does not change the resulting bivariant theory \cite{Walk:Thom} and has several advantages such as making certain induced maps of spectra functorial rather than merely functorial up to homotopy.

The projection $\mcal{O}^{n+1} \to \mcal{O}^{n}$ onto the first $n$ basis elements defines a map $G_{Y}(n) \to G_{Y}(n+1)$ and let $G_{Y}$ denote the ind-scheme
$ G_{Y} = \colim_{n} G_{Y}(n)$. The elements of $G_{Y}(X) = \colim_{n}\Hom{}{X}{G_{Y}(n)} $ are considered as quotient objects $[p:\mcal{O}^{\infty}_{X\times Y}\epi \mcal{F}]$ such that all but finitely many of the standard basis elements are in the kernel of $p$.

\begin{definition}
An $m$-tuple $([p_{1}:\mcal{O}^{\infty}_{X\times Y} \epi \mcal{F}_{1}],\ldots, [p_{m}:\mcal{O}^{\infty}_{X\times Y}\epi \mcal{F}_{m}])$ is said to be in general position if the induced map
\begin{equation*}
 \sum_{i} p_{i}: \mcal{O}^{\infty}_{X\times Y} \to \mcal{F}_{1}\oplus \cdots \oplus \mcal{F}_{m}
\end{equation*}
is again a surjection.  
Write $G_{Y}^{(m)}(X) \subseteq G_{Y}(X)^{\times m}$ for the subset 
of all $m$-tuples which are in general position.
\end{definition}
By \cite{GW:geom} the functor $X\mapsto G_{Y}^{(m)}(X)$ on $Sch/k$ is represented by an open ind-scheme $G_{Y}^{(m)}\subseteq G_{Y}^{\times m}$. By convention $G_{Y}^{(0)} = \spec k$.

The collection $\{G_{Y}^{(m)}\}_{m\geq 0}$  gives a $\Gamma$-object in presheaves of sets on $Sch/k$ (see Appendix \ref{gammaspc} for a recollection of $\Gamma$-spaces). Given a map $\sigma: \underline{m} \to \underline{n}$ then $\sigma_{*}:G_{Y}^{(m)}(X) \to G_{Y}^{(n)}(X)$ is defined by 
 \begin{align*}
 \sigma_{*}([p_{1}:\mcal{O}^{\infty}_{X\times Y} \epi \mcal{F}_{1}],& \cdots,  [p_{m}:\mcal{O}_{X\times Y}^{\infty}\epi \mcal{F}_{m}]) =\\
 = & ([q_{1}:\mcal{O}^{\infty}_{X\times Y} \epi \mcal{E}_{1}],\cdots, [q_{n}:\mcal{O}_{X\times Y}^{\infty}\epi \mcal{E}_{n}]) ,
 \end{align*}
where $\mcal{E}_{k} = \oplus_{i\in \sigma^{-1}(\underline{k})}\mcal{F}_{i}$ and $q_{k}:\mcal{O}^{\infty}_{X\times Y}\epi \mcal{E}_{k}$ is the sum over $i\in \sigma^{-1}(\underline{k})$ of the maps $p_{i}$. The condition of being in general-position implies that this definition is well-defined. This is functorial in $X$ and so defines a $\Gamma$-object in the category of presheaves on $Sch/k$. The argument in \cite[Lemma 2.2]{GW:geom} shows that the $\Gamma$-space $m\mapsto G_{Y}^{(m)}(X\times \Delta^{\bullet})$ is special for any $X$. 

\begin{definition}\label{defkgeom}
For quasi-projective $k$-varieties $X$ and $Y$ define $\mcal{K}_{geom}(X,Y)$ to be the $\Omega$-spectrum associated with the special $\Gamma$-space $m\mapsto G_{Y}^{(m)}(X\times\Delta^{\bullet})$.
\end{definition}
 By \cite{Walk:Thom} the $\Omega$-spectrum $\mcal{K}_{geom}(X,Y)$ defined here is equivalent to the spectrum obtained from the special $\Gamma$-space $m\mapsto K_{Y}^{(m)}(X\times \Delta^{\bullet})$. Let $\mcal{K}(X,Y)$ denote the $K$-theory spectrum of the exact category $\mcal{P}(X,Y)$ of coherent $\mcal{O}_{X\times Y}$-modules which are flat over $X$ and have finite support over $X$. By \cite[Theorem 2.3]{GW:geom} for any $X$ and $Y$ there is a natural map of $\Omega$-spectra $\mcal{K}_{geom}(X,Y) \to |\mcal{K}(X\times\Delta^{\bullet},Y)|$ which is always a weak equivalence. In particular whenever $X$ is regular then $\pi_{*}\mcal{K}_{geom}(X,Y) = K_{*}(X,Y)$.

\begin{lemma}(\cite[Lemma 2.3]{Walk:Thom})
A map $f: Y\to Y'$ of $k$-varieties induces a natural map of $\Gamma$-spaces $f_{*}:\mcal{G}_{Y}^{(m)}(X\times\Delta^{\bullet}) \to \mcal{G}_{Y}^{(m)}(X\times \Delta^{\bullet})$, by $f_{*}[\mcal{O}_{X\times Y}^{n}\twoheadrightarrow \mcal{M}] = [\mcal{O}_{X\times Y'}^{n} \to (id\times f)_{*}\mcal{M}]$. In particular we obtain a natural transformation of presheaves of $\Omega$-spectra
\begin{equation*}
f_{*}:\mcal{K}_{geom}( - , Y) \to \mcal{K}_{geom}(-,Y').
\end{equation*}
\end{lemma}
\begin{proof}
Let $[\mcal{O}_{X\times Y}^{n}\twoheadrightarrow \mcal{M}]$ be an object of $G_{Y}(n)(X)$. The second condition on $\mcal{M}$ in Definition \ref{Gn} guarantees that the induced map $[\mcal{O}_{X\times Y'}^{n} \to (id\times f)_{*}\mcal{M}]$ is again a surjection. Thus $f_{*}:G_{Y}(X\times \Delta^{\bullet}) \to G_{Y'}(X\times\Delta^{\bullet})$ is well-defined. The second condition in Definition \ref{Gn} also guarantees that if 
$([p_{1}:\mcal{O}_{X\times Y}^{\infty}\twoheadrightarrow \mcal{M}_{1}], \ldots, [p_{m}: \mcal{O}_{X\times Y}^{\infty}\twoheadrightarrow \mcal{M}_{m}])$ is an $m$-tuple in general position then the induced $m$-tuple $((id\times f)_{*}(p_{1}),\ldots, (id\times f)_{*}(p_{m}))$ is again in general position. Thus we have a natural map of $\Gamma$-spaces $f_{*}:G_{Y}^{(m)}(X\times\Delta^{\bullet})\to G_{Y'}^{(m)}(X\times\Delta^{\bullet})$.
\end{proof}

Let $Cube(n)$ denote the category whose objects are the subsets of $\{1,\ldots, n\}$ and morphisms are subset inclusions. An $n$-cube in $Sch/k$ is a functor from $Cube(n)$ to $Sch/k$.  
For a subset $I\subseteq \{1,\ldots, n\}$ define $\P^{\wedge n}_{I} = P_{1}\times\cdots\times P_{n}$ where $P_{i} = \P^{1}$ if $i\in I$ and otherwise $P_{i} = \{\infty\}$. 
By the above lemma we have natural transformations $\mcal{K}_{geom}(-, \P^{\wedge n}_{I}) \to \mcal{K}_{geom}(-, \P^{\wedge n}_{J})$ for $I\subseteq J$. Define $\mcal{W}^{(n)}(X)$ to be homotopy colimit of the cube $I\mapsto \mcal{K}_{geom}(-, \P^{\wedge n}_{I})$, this is natural in $X$ and so defines the presheaf $\mcal{W}^{(n)}(-)$.  When $X$ is regular, $\mcal{W}^{(0)}(X) \sim \mcal{K}(X)$.

There are functorial maps $\mcal{W}^{(n+1)}(X) \to \mcal{W}^{(n)}(X)$ \cite[Section 2]{Walker:adams}, from which we obtain a tower of presheaves
\begin{equation*}
 \cdots \to \mcal{W}^{(n+1)}(-) \to \mcal{W}^{(n)}(-) \to \cdots \to \mcal{W}^{(0)}(-).
\end{equation*}

For a simplicial abelian group $A_{\bullet}$ write $\mcal{M}(A_{\bullet})$ for the associated $\Omega$-spectrum, $\mcal{M}(A_{\bullet}) = (A_{\bullet}, BA_{\bullet}, B^{2}A_{\bullet}, \cdots)$.

Write $z_{equi}(Y,0)(-)$ for the presheaf of equidimensional cycles of relative  dimension $0$ \cite{SV:rel}. Using the localization property for the complex of equidimensional cycles \cite{FV:biv} we see that the complex $z_{equi}(\P^{\wedge n}, 0)(X\times\Delta^{\bullet})$ is equivalent to the complex $z_{equi}(\A^{n}, 0)(X\times\Delta^{\bullet})$ which computes motivic cohomology (see \cite{FV:biv,MVW:mot}). As explained in \cite[Section 4]{Walker:adams} there is a natural sequence of spectra
\begin{equation*}
 \mcal{W}^{(n+1)}(X) \to \mcal{W}^{(n)}(X) \to \mcal{M}(z_{equi}(\P^{\wedge n}, 0)(X\times\Delta^{\bullet}))
\end{equation*}
which is homotopic to the constant map.  By \cite{Suslin:GraySS} this sequence is a homotopy fiber sequence when $X = \spec \mcal{O}_{Y,y}$ is the local ring of a smooth variety $Y$ at a point $y\in Y$.

This tower is globalized as in Friedlander-Suslin \cite{FS:AHSS} via Godement resolutions (for the Zariski topology) which are recalled in Appendix \ref{glob}. Write $\mcal{G}$ for the Godement resolution and define $\mcal{K}^{(n)}(-) = \mcal{G}\mcal{W}^{(n)}(-)$ and  $\mcal{M}^{(n)}(-) = \mcal{G}\mcal{M}(z_{equi}(\P^{\wedge n}, 0)(-\times\Delta^{\bullet}))$. 
Observe that because the presheaf of simplicial abelian groups $z_{equi}(\P^{\wedge n},0)(-\times\Delta^{\bullet})$ satisfies Zariski descent on smooth varieties, the natural map $\mcal{M}(z_{equi}(\P^{\wedge n}, 0)(X\times\Delta^{\bullet})) \to \mcal{M}^{(n)}(X)$ is a weak equivalence for any smooth $k$-variety $X$.

Because $\mcal{W}^{(n+1)}(X) \to \mcal{W}^{(n)}(X) \to \mcal{M}(z_{equi}(\P^{\wedge n}, 0)(X\times\Delta^{\bullet}))$ is a homotopy fiber sequence when $X = \spec\mcal{O}_{Y,y}$ for a smooth variety $Y$ and a point $y\in Y$, we see that 
\begin{equation}\label{hofibss}
\mcal{K}^{(n+1)}(X) \to \mcal{K}^{(n)}(X) \to \mcal{M}^{(n)}(X)
\end{equation}
is a homotopy fiber sequence for any smooth $X$. In general when $X$ is not smooth the composition is homotopic to the constant map.

This gives the tower of presheaves of spectra on $Sch/k$
\begin{equation}\label{mottow}
 \xymatrix@-1pc{
\cdots \ar[r] & \mcal{K}^{(q+1)} \ar[r] & \mcal{K}^{(q)} \ar[r]\ar[d] & \mcal{K}^{(q-1)} \ar[r]\ar[d] & \cdots \ar[r] & \mcal{K}^{(0)} \ar[d]\ar[r] & \mcal{K} . \\
 & & \mcal{M}^{(q)} & \mcal{M}^{(q-1)} & & \mcal{M}^{(0)} &
}
\end{equation}

We summarize the properties of this tower which we use.

\begin{theorem}\label{mtow}\cite{Gray:wt,Suslin:GraySS}
Let $X$ be a smooth  $k$-variety.
\begin{enumerate}
\item Each of $\mcal{K}^{(q)}(X)$ and $\mcal{M}^{(q)}(X)$ are  $(-1)$-connected $\Omega$-spectra.

\item The sequence
\begin{equation*}
 \mcal{K}^{(q+1)}(X) \to \mcal{K}^{(q)}(X) \to \mcal{M}^{(q)}(X)
\end{equation*}
is a homotopy fibration sequence of spectra (in general the sequence is homotopic to the constant map if $X$ is not smooth).

\item For any abelian group $A$,
\begin{equation*}
 \pi_{k}(\mcal{M}^{(q)}(X); A) = H_{\mcal{M}}^{2q-k}(X; A(q)) 
\end{equation*}

\item The augmentation $\mcal{K}^{(0)}(X) \to \mcal{K}(X)$ is a weak equivalences of spectra.

\item If $k < q-\dim X$ then $\pi_k \mcal{K}^{(q)}(X) = 0$.

\end{enumerate}
\end{theorem}
\begin{proof}
The first four items have already been discussed above. Let $\widetilde{\pi}_{i}\mcal{W}^{(q)}$ denote the Zariski sheafification of the presheaf $U\mapsto \pi_{i}\mcal{W}^{(q)}(U)$ on $Sm/k$. For $Y=\spec\mcal{O}_{X,x}$, the localization of a smooth variety $X$ a point $x\in X$, there is a natural homotopy equivalence of spectra $\mcal{K}_{geom}(Y,\mathbb{G}_{m}^{\wedge q}) \to \Omega^{q}\mcal{W}^{(q)}(Y)$ by \cite[Theorem 7.11]{Walk:thes}. Therefore $\widetilde{\pi}_{i}\mcal{W}^{(q)} = 0$ for $i < q$. That $\pi_k \mcal{K}^{(q)}(X) = 0$ for $k< q- \dim X$ now follows from the descent spectral sequence \cite{BG} 
$$
E_{2}^{s,t}= H^{s}_{Zar}(X;\widetilde{\pi}_{t}\mcal{W}^{(q)}) \Longrightarrow \pi_{t-s}\mcal{K}^{(q)}(X)
$$
since $E_{2}^{s,t} \neq 0$ means that $s\leq \dim X$ and $t\geq q$.
\end{proof}

\begin{corollary}\cite{Gray:wt,Suslin:GraySS}
For any smooth  $k$-variety $X$ and any abelian group $A$, there is a strongly convergent spectral sequence
\begin{equation*}
E_{2}^{p,q}= H^{p-q}_{\mcal{M}}(X ; A(-q)) \Longrightarrow K_{-p-q}(X; A).
\end{equation*} 
\end{corollary}
\begin{proof}
For any abelian group $A$, the homotopy fiber sequences (\ref{hofibss}) yield an exact couple,
\begin{equation*}
 D_{2}^{p,q} = \pi_{-p-q}(\mcal{K}^{(-q)}(X) ; A) \;\;\; \text{and} \;\;\; E_{2}^{p,q} = \pi_{-p-q}(\mcal{M}^{(-q)}(X) ; A).
\end{equation*}
By the previous theorem this exact couple is bounded below and thus the resulting spectral sequence converges to $\colim_{p}D_{2}^{p,n-p} = \pi_{-n}(\mcal{K}^{(p-n)}(X) ; A) = \pi_{-n}\mcal{K}(X;A)$ (see for example \cite[Corollary 5.9.7]{Weibel:hom}).
\end{proof}

We are primarily interested in the cases $k=\R$ and $k=\C$. Let $L/k$ be a finite extension of fields. In the following lemma, in order to (temporarily) distinguish between the presheaf of spectra on $Sch/k$ obtained using $G_{\P^{n}_{k}}$ and the presheaf on $Sch/L$ using $G_{\P^{n}_{L}}$ we use the notation $\mcal{K}_{k}^{(n)}(X)$ and $\mcal{K}_{L}^{(n)}(X)$ (and similarly for the presheaves $\mcal{M}^{(n)}$). Because of the following lemma we will not make careful distinction between the two in the rest of this paper.
 
\begin{lemma}\label{bch}
Let $L/k$ be a finite field extension and let $p:Sch/k \to Sch/L$ be the map of sites specified by $X\mapsto X_{L}$. Then for an $L$-variety $Y$ we have that $p^{*}\mcal{K}_{k}^{(n)}(Y)=\mcal{K}^{(n)}_{L}(Y)$ and $p^{*}\mcal{M}_{k}^{(n)}(Y)=\mcal{M}^{(n)}_{L}(Y)$ 
\end{lemma}
\begin{proof}
For a $k$-variety $Y$ we have the agreement of presheaves of sets on $Sch/L$,
$$
p^{*}\Hom{k}{-}{G_{Y}^{(m)}} = \Hom{L}{-}{G_{Y_{L}}^{(m)}}.
$$ 
Thus for any $L$-variety $X$ we have, $p^{*}G_{Y}^{(m)}(X\times_{L}\Delta_{L}^{\bullet}) = G_{Y_{L}}^{(m)}(X\times_{L}\Delta_{L}^{\bullet})$. If $F(-)$ is a presheaf of $\Gamma$-spaces then $p^{*}BF = Bp^{*}F$ and so  $p^{*}\mcal{K}_{geom, k}(-,Y) = \mcal{K}_{geom, L}(-, Y_{L})$. This implies that
\begin{equation*}
p^{*}\mcal{W}^{(n)}_{k}(-) = \mcal{W}^{(n)}_{L}(-).
\end{equation*} 
Finally, because the Godement resolution commutes with $p^{*}$ we have
\begin{equation*}
p^{*}\mcal{K}^{(n)}_{k}(-) = p^{*}\mcal{G}\mcal{W}^{(n)}_{k}(-) = \mcal{G}\mcal{W}^{(n)}_{L}(-) = \mcal{K}^{(n)}_{L}(-).
\end{equation*}
A similar argument gives the result for the presheaves $\mcal{M}^{(n)}(-)$.
%
%
\end{proof}

\subsection{Semi-topological Spectral Sequence}\label{sstss}
We now construct a semi-topological spectral sequence for real varieties, using the motivic spectral sequence above, and following the method of \cite{FHW:sst}. We begin by recalling a construction used by Friedlander-Walker.

Let $F(-)$ be a presheaf of abelian groups, simplicial sets, or spectra on $Sch/\R$. If $T$ is a topological space and $X$ is a real variety then define $F(T\times_{\R}X)$ by the filtered colimit
\begin{equation*}
 F(T\times_{\R}X) = \colim_{T\to V(\R)} F(V\times_{\R} X).
\end{equation*}
The colimit is over continuous maps $f:T\to V(\R)^{an}$, where $V$ is a real variety, and a map from $(f:T\to V(\R))$ to $(g:T\to U(\R))$ is a map $\alpha:V\to U$ of real varieties such that $g=\alpha f$. 
Similarly define $F(T\times_{\C} X)$ by the colimit
\begin{equation*}
 F(T\times_{\C}X) = \colim_{T\to V(\R)} F((V_{\C})\times_{\R} X).
\end{equation*}
Via the action of $Gal(\C/\R)$ on each $F((V_{\C})\times_{\R} X)$ we obtain an action of $Gal(\C/\R)$ on $F(T\times_{\C}X)$.

\begin{remark}\label{base}
When $H$ is a presheaf on $Sch/\C$, Friedlander-Walker \cite[proof of Lemma 2.4]{FW:real} show that 
$$
\colim_{T\to U(\R)}H(U_{\C})\iso \colim_{T\to W(\C)}H(W)
$$
where the first colimit is over continuous maps $T\to U(\R)$ with $U$ a real variety and in the second $T\to W(\C)$ with $W$ a complex variety. In particular writing $p:Sch/\R\to Sch/\C$ for the map $X\mapsto X_{\C}$ we see that for a presheaf $F(-)$ on $Sch/\R$ we have the agreement 
 $$
 F(T\times_{\C}\R) = (p_{*}p^{*}F)(T\times_{\R} \R) = (p^{*}F)(T\times_{\C}\C).
 $$ 
\end{remark}

Let $F$ be a presheaf of simplicial sets on $Sch/\R$. We obtain a bisimplicial set  by the formula $ d\mapsto F(\Delta^{d}_{top}\times_{\R} X)$
and define
\begin{equation*}
 F(\Delta^{\bullet}_{top}\times_{\R} X) = \diag(d\mapsto F(\Delta^{d}_{top}\times_{\R} X)) .
\end{equation*}

Similarly, if $\mcal{F}$ is a presheaf of spectra on $Sch/\R$ then $d\mapsto F(\Delta^{d}_{top}\times_{\R} X)$ is a simplicial spectrum. Define the spectrum $F(\Delta^{\bullet}_{top}\times_{\R} X)$ by the formula 
$$
F(\Delta^{\bullet}_{top}\times_{\R}X)^{i} = \diag(d\mapsto F(\Delta^{d}_{top}\times_{\R} X)^{i})
$$ 
and bonding maps given by
\begin{align*}
 S^{1}\wedge \diag(d\mapsto F(\Delta^{d}_{top}\times_{\R}X)^{i})  \iso & \diag(d\mapsto S^{1}\wedge F(\Delta^{d}_{top}\times_{\R} X)^{i}) \to \\
 \to & \diag(d\mapsto F(\Delta^{d}_{top}\times_{\R}X)^{i+1})
\end{align*}

\begin{definition}
Let $F$ be a presheaf of simplicial sets of spectra. Define a new presheaf  $F_{sst}$ on $Sch/\R$ by
\begin{equation*}
 F_{sst}(X) = F(\Delta^{\bullet}_{top}\times_{\R} X).
\end{equation*}
\end{definition}

\begin{remark}
For emphasis and to avoid confusion we will often write below $F(\Delta^{\bullet}_{top}\times_{\R}\R)$ instead of $F(\Delta^{\bullet}_{top})$ when $F$ is a presheaf on $Sch/\R$ and similarly we will write $G(\Delta^{\bullet}_{top}\times_{\C}\C)$ instead of $G(\Delta^{\bullet}_{top})$ when $G$ is a presheaf on $Sch/\C$.
\end{remark}

In \cite{FW:ratisos} the $uad$-topology on $Sch/\C$ is used to establish a recognition principle identifying when a map $F(-)\to G(-)$ of presheaves on $Sch/\C$ induces a weak equivalence $F(\Delta^{\bullet}_{top}\times_{\C}\C)\to G(\Delta^{\bullet}_{top}\times_{\C}\C)$. In \cite{HV:VT} it is observed that their work carries over to $Sch/\R$, using a suitable real $uad$-topology. We briefly recall this topology, which is essentially due to Deligne.
\begin{defn}
\begin{enumerate}
\item A continuous map of topological spaces $f: S\to T$ is said to satisfy \textit{cohomological descent} if for any sheaf $A$ of abelian groups on $T$ the natural map 
\begin{equation*}
 H^*(T, A) \to H^{*}(N_T(S), f^*A)
\end{equation*}
is an isomorphism. Here $N_T(S)\to T$ is the Cech nerve of $f$, i.e. $N_T(S)$ is the simplicial space which in degree $n$ is the $n+1$-fold fiber product of $S$ over $T$.
A map $f:S\to T$ is said to be of \textit{universal cohomological descent} provided the pullback $S\times_{T}T'\to T'$ along any continuous map $T'\to T$ is again of cohomological descent.

\item The $uad$-topology on $Sch/{\C}$ is the Grothendieck topology associated to the pretopology generated by collections $\{U_i\to X\}$ such that $\coprod U_i(\C)^{an}\to X(\C)^{an}$ is a surjective map of universal cohomological descent.

\item The $uad$-topology on $Sch/{\R}$ is the Grothendieck topology associated to the pretopology generated by collections $\{U_i\to X\}$ such that $\coprod U_i(\R)^{an}\to X(\R)^{an}$ is a surjective map of universal cohomological descent.
\end{enumerate}
\end{defn}

A useful feature of the $uad$-topologies is that a complex variety is locally smooth in the $uad$-topology and a real variety is locally smooth in the real $uad$-topology.  

The following is the recognition principle proved by Friedlander-Walker \cite{FW:ratisos} (for the real analogue see \cite{HV:VT}).
\begin{theorem}[{\cite[Theorem 2.2]{FW:ratisos}}]\label{recog}
Let $k=\C$ or $\R$. Suppose that $F\to G$ is a natural transformation of presheaves of abelian groups on $Sch/{k}$. If $F_{uad} \xrightarrow{} G_{uad}$ is an isomorphism of $uad$-sheaves, then 
\begin{equation*}
 F(\Delta^{\bullet}_{top}\times_{k} k) \to G(\Delta^{\bullet}_{top}\times_{k} k)
\end{equation*}
 is a homotopy equivalence of simplicial abelian groups.
\end{theorem}

The following two propositions were proved in the complex case in \cite{FHW:sst}. 
\begin{proposition}(c.f. \cite[Theorem 2.6]{FHW:sst}\label{srec}
Let $k= \C$ or $\R$. Suppose that $f:\mcal{F}\to \mcal{G}$ is a natural transformation of presheaves of spectra on  $Sch/k$. Suppose that there is an integer $N$ with the property that  the maps $(\pi_{q}\mcal{F})_{uad} \to (\pi_{q}\mcal{G})_{uad}$ are isomorphisms of $uad$-sheaves for $q\leq N$ and a surjection for $q=N+1$. Then 
\begin{equation*}
\pi_{q} \mcal{F}(\Delta^{\bullet}_{top}\times_{k} k) \to \pi_{q}\mcal{G}(\Delta^{\bullet}_{top}\times_{k} k)
\end{equation*}
 is an isomorphism for $q\leq N$ and a surjection for $q=N+1$ . In particular if $\mcal{F}(X)\to \mcal{G}(X)$ is a weak equivalence of spectra for all smooth $X$ then 
 $$
 \mcal{F}(\Delta^{\bullet}_{top}\times_{k} k) \to \mcal{G}(\Delta^{\bullet}_{top}\times_{k} k)
 $$ 
 is also a weak equivalence of spectra.
\end{proposition}
\begin{proof}
Define the presheaf of spectra $\mcal{Q}(-)$ by $\mcal{Q}(X) = \hocofib(\mcal{F}(X) \to \mcal{G}(X))$. Then $\pi_{q}\mcal{Q}_{uad} = 0$ for all $q\leq N+1$. By the recognition principle the simplicial abelian group  $(d\mapsto \pi_{q}\mcal{Q}(\Delta^{d}_{top}\times_{k} k))$  is contractible for any $q\leq N+1$.

From the strongly convergent spectral sequence \cite[Corollary 4.22]{Jar:et}
\begin{equation*}
 E^{2}_{s,t} = \pi_{s}(d\mapsto\pi_{t}\mcal{Q}(\Delta^{d}_{top}\times_{k} k)) \Longrightarrow \pi_{s+t}\mcal{Q}(\Delta^{\bullet}_{top}\times_{k} k)
\end{equation*}
we conclude that $\pi_{q}\mcal{Q}(\Delta^{\bullet}_{top}\times_{k} k) = 0$ for $q\leq N+1$. Finally observe that 
$$
\mcal{Q}(\Delta^{d}_{top}\times_{k} k) \wkeq \hocofib(\mcal{F}(\Delta^{d}_{top}\times_{k} k)\to\mcal{G}(\Delta^{d}_{top}\times_{k} k))
$$
because filtered colimits preserve homotopy cofiber sequences. Since the diagonal converts degreewise homotopy cofiber sequences of simplicial spectra into homotopy cofiber sequences of spectra, we have that 
$$
\mcal{Q}(\Delta^{\bullet}_{top}\times_{k} k) \simeq \hocofib(\mcal{F}(\Delta^{\bullet}_{top}\times_{k} k) \to \mcal{G}(\Delta^{\bullet}_{top}\times_{k} k)).
$$ 
Finally from the long exact sequence for homotopy cofiber sequences of spectra we conclude that 
$\pi_{k} \mcal{F}(\Delta^{\bullet}_{top}\times_{k} k) \to \pi_{k}\mcal{G}(\Delta^{\bullet}_{top}\times_{k} k)$
 is an isomorphism for $k\leq N$ and a surjection for $k=N+1$.
\end{proof}

\begin{proposition}(\cite[Corollary 2.7]{FHW:sst})\label{hfseq}
 Let $k=\C$ or $\R$. Suppose that
\begin{equation*}
\mcal{F} \to \mcal{G} \to \mcal{H}
\end{equation*}
is a sequence of presheaves of spectra on $Sch/k$ such that the composition is homotopic to the constant map and the sequence 
\begin{equation*}
\mcal{F}(X) \to \mcal{G}(X) \to \mcal{H}(X)
\end{equation*}
is a homotopy fiber sequence for any smooth $X$. Then
\begin{equation}\label{cof1}
\mcal{F}(\Delta^{\bullet}_{top}\times_{k} X) \to \mcal{G}(\Delta_{top}^{\bullet}\times_{k} X) \to \mcal{H}(\Delta^{\bullet}_{top}\times_{k} X)
\end{equation}
is a homotopy fiber sequence of spectra.
\end{proposition}
\begin{proof}
A sequence of spectra is a homotopy fiber sequence if and only if it is a homotopy cofiber 
sequence. Write $\mcal{Q}(X) = \hocofib(\mcal{F}(X) \to \mcal{G}(X))$.  Filtered colimits preserve homotopy cofiber sequences of spectra. Therefore for each $d$ the 
sequence
\begin{equation*}
 \mcal{F}(\Delta^{d}_{top}\times_{k} k) \to \mcal{G}(\Delta_{top}^{d}\times_{k} k) \to 
\mcal{Q}(\Delta^{d}_{top}\times_{k} k)
\end{equation*}
is a homotopy cofiber sequence. The diagonal converts degreewise homotopy cofiber sequences 
of simplicial spectra into homotopy cofiber sequences of spectra and thus
\begin{equation*}
 \mcal{F}(\Delta^{\bullet}_{top}\times_{k} k) \to \mcal{G}(\Delta_{top}^{\bullet}\times_{k} k) \to 
\mcal{Q}(\Delta^{\bullet}_{top}\times_{k} k)
\end{equation*}
is a homotopy cofiber sequence. The natural map 
$\mcal{Q}\to\mcal{H}$ is a weak equivalence on smooth varieties which implies that $(\pi_{q}\mcal{Q})_{uad}\to(\pi_{q}\mcal{H})_{uad}$ is an isomorphism for all $q$. Therefore 
$\mcal{Q}(\Delta^{\bullet}_{top}\times_{k} k) \to \mcal{H}(\Delta^{\bullet}_{top}\times_{k} k)$ is also a weak equivalence by Proposition \ref{srec}. We conclude that the sequence (\ref{cof1}) is a homotopy cofiber (and thus a homotopy fiber) sequence of spectra.
\end{proof}

Write $\mcal{K}\R^{(q)}_{sst}(X) = \mcal{K}^{(q)}(X\times_{\R}\Delta^{\bullet}_{top})$ and $\mcal{M}\R^{(q)}_{sst}(X)=\mcal{M}^{(q)}(X\times_{\R}\Delta^{\bullet}_{top})$. We obtain from the motivic tower (\ref{mottow}) the semi-topological tower of presheaves of spectra,
\begin{equation*}
 \xymatrix@-1pc{
\cdots \ar[r] & \mcal{K}\R_{sst}^{(q+1)} \ar[r] & \mcal{K}\R_{sst}^{(q)} \ar[r]\ar[d] & \mcal{K}\R_{sst}^{(q-1)} \ar[r]\ar[d] & \cdots \ar[r] & \mcal{K}\R_{sst}^{(0)} \ar[d]\ar[r] & \mcal{K}\R_{sst} . \\
 & & \mcal{M}\R_{sst}^{(q)} & \mcal{M}\R_{sst}^{(q-1)} & & \mcal{M}\R_{sst}^{(0)} &
}
\end{equation*}

\begin{theorem}\label{sstowtprop}
Let $X$ be a smooth real variety.
\begin{enumerate}

\item The spectra $\mcal{K}\R_{sst}^{(q)}(X)$ and $\mcal{M}\R_{sst}^{(q)}(X)$ are $(-1)$-connected spectra. Moreover  $\pi_k \mcal{K}\R_{sst}^{(q)}(X) = 0$ for $k< q- \dim X$.

\item The sequence
\begin{equation}\label{ssthofib}
 \mcal{K}\R_{sst}^{(q+1)}(X) \to \mcal{K}\R_{sst}^{(q)}(X) \to \mcal{M}\R_{sst}^{(q)}(X)
\end{equation}
is a homotopy fiber sequence of spectra.

\item For any abelian group $A$,
\begin{equation*}
 \pi_{k}(\mcal{M}\R_{sst}^{(q)}(X); A) = L^{q}H\R^{q-k, q}(X; A) ,
\end{equation*}
where $L^{q}H\R^{q-k,q}(X;A)$ denotes the real morphic cohomology defined in \cite{FW:real}.
\item The augmentation
\begin{equation*}
 \mcal{K}\R_{sst}^{(0)}(X) \to \mcal{K}\R_{sst}(X)
\end{equation*}
is a weak equivalence of spectra, where $\mcal{K}\R_{sst}(X)$ is the real semi-topological $K$-theory spectrum defined in \cite[Definition 2.1]{FW:real} (denote there by $\mcal{K}^{alg}(\Delta^{\bullet}_{top}\times_{\R}X)$.
\end{enumerate}
\end{theorem}
\begin{proof}
 \begin{enumerate}
  \item This follows from Theorem \ref{mtow} and Proposition \ref{srec}.
\item This follows from Theorem \ref{mtow} and Proposition \ref{hfseq}.
\item The spectrum $\mcal{M}^{(q)}(X)$ is the spectrum associated to the simplicial abelian group $z_{equi}(\P_{\R}^{\wedge q},0)(X_{\R}\times_{\R}\Delta_{\R}^{\bullet})$. Consequently we have that $\mcal{M}\R^{(q)}_{sst}(X)$ is naturally equivalent to the spectrum associated to the simplicial abelian group $z_{equi}(\A_{\R}^{ q},0)(X_{\R}\times\Delta_{\R}^{\bullet}\times\Delta^{\bullet}_{top})$. By \cite{HV:VT} the homotopy groups of this simplicial abelian group compute real morphic cohomology. 
\item Follows from Theorem \ref{recog} and Theorem \ref{mtow}.
 \end{enumerate}
\end{proof}

The following is the Atiyah-Hirzebruch type spectral sequence relating the real morphic cohomology and the real semi-topological algebraic $K$-theory, of a smooth real variety $X$.
\begin{corollary}
For any smooth $k$-variety $X$ and any abelian group $A$, there is a strongly convergent spectral sequence
\begin{equation*}
E_{2}^{p,q}= L^{-q}H\R^{p,-q}(X ; A) \Longrightarrow K\R^{sst}_{-p-q}(X; A).
\end{equation*} 
\end{corollary}
\begin{proof}
For any abelian group $A$, the homotopy fiber sequences (\ref{ssthofib}) yield an exact couple,
\begin{equation*}
 D_{2}^{p,q}(sst) = \pi_{-p-q}(\mcal{K}\R_{sst}^{(-q)}(X) ; A) \;\;\; \text{and} \;\;\; E_{2}^{p,q}(sst) = \pi_{-p-q}(\mcal{M}\R_{sst}^{(-q)}(X) ; A).
\end{equation*}
This exact couple is bounded below by the previous theorem and therefore the resulting spectral sequence converges to $\colim_{p}D_{2}^{p,n-p} =  \pi_{-n}\mcal{K}\R_{sst}(X;A)$ (see for example \cite[Corollary 5.9.7]{Weibel:hom}).
\end{proof}

\begin{remark}
More generally we could use the equivariant homotopy groups $\pi_{s,t}$ rather than $\pi_{s,0}$ in order to produce other spectral sequences relating equivariant morphic cohomology and equivariant semi-topological $K$-theory $\mcal{K}\R_{s,t}(X) \eqdef \pi_{s,t}\mcal{K}(X_{\C}\times_{\C} \Delta^{\bullet}_{top})$. 
\end{remark}

\subsection{Topological spectral sequence}\label{topss}
We continue to write $G=Gal(\C/\R)$. Let $F$ be a presheaf of simplicial sets or spectra on $Sch/\R$. We remind the reader that for a topological space $T$ we have a natural $G$-action on $F(T\times_{\C}X) = \colim_{T\to U(\R)}F(U_{\C}\times_{\R} X)$
and that $F(T\times_{\C} X) = p^{*}F(T\times_{\C} X_{\C})$ where $p:Sch/\R \to Sch/\C$ is the map of sites given by $X\mapsto X_{\C}$ (see Remark \ref{base}).

The construction of \cite[Lemma 3.2]{FHW:sst} yields an equivariant natural transformation of presheaves of $G$-simplicial sets
\begin{equation*}
 p^{*}F(\Delta^{\bullet}_{top}\times_{\C} X_{\C}) \to \shom{}{\sing_{\bullet}X(\C)}{p^{*}F(\Delta^{\bullet}_{top})},
\end{equation*}
where $\shom{}{-}{-}$ denotes the simplicial function complex equipped with the usual $G$-action. If $F$ is a presheaf of spectra then this is an equivariant natural transformation of naive $G$-spectra.

The projections $U_{\C} \to U$ induce the natural map
\begin{equation*}
F(T\times_{\R}X)=\colim_{T\to U(\R)}F(U\times_{\R} X)\xrightarrow{}\left(\colim_{T\to U(\R)}F(U_{\C}\times_{\R} X)\right)^{G}= p^{*}F(T\times_{\C} X_{\C})^{G} . 
\end{equation*}
Composing with the map from \cite[Lemma 3.2]{FHW:sst} yields the natural transformation
\begin{equation*}
F(\Delta_{top}^{\bullet}\times_{\R} X) \to \shom{G}{\sing_{\bullet}X(\C)}{p^{*}F(\Delta^{\bullet}_{top})},
\end{equation*}
where $\shom{G}{-}{-}$ is the simplicial complex of equivariant maps. 

\begin{definition}
If $F$ is a presheaf of simplicial sets on $Sch/\R$ define a presheaf of $G$-simplicial sets on $Sch/\R$ by
\begin{equation*}
F_{top}(X) = \shom{}{\sing_{\bullet}X(\C)}{p^{*}F(\Delta^{\bullet}_{top})}.
\end{equation*}
If $F$ is a presheaf of spectra on $Sch/\R$ define $F_{top}$ by the same formula. In this case $F_{top}$ is a presheaf of naive $G$-spectra.
\end{definition}
The preceding discussion shows that we have the sequence of natural transformations 
\begin{equation*}
 F \to F_{sst} \to F_{top}^{G}
\end{equation*}
of presheaves  on $Sch/\R$.

Applying the functor $(-)_{top}$ to the motivic tower on $Sch/\R$ gives rise to the tower of presheaves of naive $G$-spectra on $Sch/\R$

\begin{equation*}
 \xymatrix@-1pc{
\cdots \ar[r] & \mcal{K}_{top}^{(q+1)} \ar[r] & \mcal{K}_{top}^{(q)} \ar[r]\ar[d] & \mcal{K}_{top}^{(q-1)} \ar[r]\ar[d] & \cdots \ar[r] & \mcal{K}_{top}^{(0)} \ar[d]\ar[r] & \mcal{K}_{top} . \\
 & & \mcal{M}_{top}^{(q)} & \mcal{M}_{top}^{(q-1)} & & \mcal{M}_{top}^{(0)} &
}
\end{equation*} 

We have taken the Grayson spectral sequence as our model for the motivic spectral sequence  because with this choice it is straightforward to see that $\mcal{K}^{(n)}(X_{\C})^{G} = \mcal{K}^{(n)}(X)$. This observation is crucial in the following proposition, which together with the two propositions after it, form the basis of the companion topological spectral sequence to the previously constructed semi-topological spectral sequence.
\begin{proposition}
Let $X$ be a smooth real variety. Then
\begin{equation}\label{tophofib}
\mcal{K}_{top}^{(q+1)}(X) \to \mcal{K}_{top}^{(q)}(X) \to \mcal{M}_{top}^{(q)}(X) 
\end{equation}
is an equivariant homotopy fiber sequence.
\end{proposition}
\begin{proof}
It is enough to verify the claim for $X = \spec\R$ and to verify that the sequence is a homotopy fiber sequence upon forgetting the $G$-action and upon taking fixed points. Evaluating this sequence on $\spec\R$ and using the observation at the end of Remark \ref{base} and Lemma \ref{bch}, the sequence becomes 
$$
\mcal{K}^{(q+1)}(\Delta^{\bullet}_{top}\times_{\C}\C) \to \mcal{K}^{(q)}(\Delta^{\bullet}_{top}\times_{\C}\C) \to \mcal{M}^{(q)}(\Delta_{top}^{\bullet}\times_{\C}\C).
$$
By Proposition \ref{gres}, upon taking fixed points this sequence becomes 
$$
\mcal{K}^{(q+1)}(\Delta^{\bullet}_{top}\times_{\R}\R) \to \mcal{K}^{(q)}(\Delta^{\bullet}_{top}\times_{\R}\R) \to \mcal{M}^{(q)}(\Delta_{top}^{\bullet}\times_{\R}\R).
$$
That both of these are homotopy fiber sequences of spectra now follows from Proposition \ref{hfseq} and Theorem \ref{mtow}.
\end{proof}

Recall that Atiyah's real $K$-theory $KR(W)$ of the space $W$ equipped with involution \cite{A:kreal} is represented by the $G$-space $\Z\times BU$ (where $BU$ has $G$-action given by complex conjugation). For a special $\Gamma$-space $F$ we write $\mathbb{S}F$ for the associated $\Omega$-spectrum. If $F$ has a $G$-action compatible with the $\Gamma$-space structure then $\mathbb{S}F$ is a naive $G$-spectrum (see Appendix \ref{gammaspc} for a recollection of these matters). We define the naive $G$-spectrum spectrum $kr$ to be the the spectrum $kr =  \mathbb{S}Grass_{\C}(\C)^{an}$ (where $Grass_{\C}(\C)^{an}$ is given the structure of a $\Gamma$-space in a manner analogous to the $\Gamma$-space structure on $G_{Y}(-\times\Delta^{\bullet})$ in Section \ref{motss}). For a $G$-space $W$ we write $kr^{i}(W) = kr^{i,0}$. Note that for $n\geq 0$ we have that $kr^{-n}(W) = \pi_{n}\mapG{W}{kr_{0}}$ where $\mapG{-}{-}$ denotes the space of equivariant maps. Furthermore we have that  $\Z\times BU = kr_{0}$ and therefore $kr^{-n}(W) = KR^{-n}(W)$ for $n\geq 0$. 
 
\begin{proposition}\label{ktopkr}
The naive $G$-spectrum $\mcal{K}_{top}(\R)$ is equivariantly equivalent to $kr$. In particular $\mcal{K}_{top}(\R)_{0} \wkeq \Z\times BU$. 
\end{proposition}
\begin{proof}
The map 
\begin{equation*}
\mcal{K}^{(0)}_{top}(\R) = \mathbb{S}\Hom{\C}{\Delta^{\bullet}_{top}\times_{\C}\Delta^{\bullet}_{\C}}{Grass_{\C}} \xleftarrow{}\mathbb{S}\Hom{\C}{\Delta^{\bullet}_{top}}{Grass_{\C}}
\end{equation*}
is a homotopy equivalence by \cite[Lemma 1.2]{FW:compK}. On fixed points 
\begin{equation*}
\mcal{K}^{(0)}_{top}(\R)^{G} = \mathbb{S}\Hom{\R}{\Delta^{\bullet}_{top}\times_{\R}\Delta^{\bullet}_{\R}}{Grass_{\R}} \xleftarrow{}\mathbb{S}\Hom{\R}{\Delta^{\bullet}_{top}}{Grass_{\R}}
\end{equation*}
is a homotopy equivalence by \cite{FW:real}.  Thus $\mathbb{S}\Hom{\C}{\Delta^{\bullet}_{top}}{Grass_{\C}}\to \mcal{K}^{(0)}_{top}(\R) $ is an equivariant equivalence. In general if $T$ is a compact Hausdorff topological space and $Y$ is a real variety then 
$\Hom{\C}{T}{Y_{\C}}\iso \Hom{cts}{T}{Y(\C)}$ is an equivariant isomorphism (by an argument similar to \cite[Proposition 4.2]{FW:sstfct}). Therefore  
$\Hom{\C}{\Delta^{\bullet}_{top}}{Grass^{(n)}_{\C}}\iso \sing_{\bullet}Grass^{(n)}_{\C}(\C)^{an}$ is an isomorphism of equivariant $\Gamma$-spaces and so $\mcal{K}^{(0)}_{top}(\R)$ is equivariantly equivalent to the the $\Omega$-spectrum $\mathbb{S}Grass_{\C}(\C)^{an}= kr$.
\end{proof}

\begin{proposition}\label{mtopem}
The $0$th space of the $\Omega$-spectrum $\mcal{M}^{(q)}_{top}(\R)$ is naturally weakly equivalent to $\sing_{\bullet}\mcal{Z}_{0}(\A^{q}_{\C})$ for $g\geq 0$. In particular for any real variety $X$,
\begin{equation*}
 \pi_{k}\mcal{M}_{top}^{(q)}(X) = H^{q-k, q}(X(\C) ; \underline{\Z}).
\end{equation*}
\end{proposition}
\begin{proof}
The spectrum $\mcal{M}(\R)^{(q)}_{top}$ is equivalent to the spectrum associated to the simplicial $G$-module $z_{equi}(\A^{q}_{\C},0)(\Delta^{\bullet}_{top})$. By \cite[Corollary 3.5]{FW:ratisos} and \cite[Proposition 4.15]{HV:VT} this is equivariantly homotopic to $\sing_{\bullet}\mcal{Z}_{0}(\A^{q}_{\C})$. By \cite{DS:equiDT} the space $\mcal{Z}_{0}(\A^{q}_{\C})$ is the Eilenberg-MacLane space $K(\underline{\Z},\R^{q,q})$.
\end{proof}

That this spectral sequence associated to the fiber sequences (\ref{tophofib}) is strongly convergent spectral sequence will follows from the following proposition.
\begin{proposition}\label{bvan}
Let $W$ be a finite-dimensional $G$-$CW$ complex of dimension $d$. There are natural isomorphisms, for any abelian group $A$
\begin{equation*}
 H^{n,m}(W;\underline{A}) \to H^{n,m+1}(W;\underline{A})
\end{equation*}
for $n+m > d$ and for $n+m+1 <0$. In particular
\begin{enumerate}
\item if $m\geq 0$, then $H^{n,m}(X,\underline{A}) = 0$ for $n> d$ and for $n+1< -m$,
\item if $m\leq 0$ then $H^{n,m}(X,\underline{A}) = 0$ for $n+m> d$ and for $n< 0$.
\end{enumerate}
\end{proposition}
\begin{proof}
Consider the homotopy cofibration sequence
\begin{equation*}
 W_{+}\wedge \Z/2_{+} \to W_{+} \to W_{+}\wedge S^{0,1} .
\end{equation*}
This induces a long exact sequence in reduced Bredon cohomology
\begin{equation*}
 \cdots \to \tilde{H}^{n-1,m+1}(W_+\wedge \Z/2_{+}) \to \tilde{H}^{n,m+1}(W_{+}\wedge S^{0,1}) \to \tilde{H}^{n,m+1}(W_{+}) \to \cdots,
\end{equation*}
which may be rewritten as
\begin{equation}\label{forles}
 \to H^{(n-1)+(m+1)}_{sing}(W;A) \to H^{n,m}(W;\underline{A}) \to H^{n,m+1}(W;\underline{A}) \to H_{sing}^{n+m+1}(W;A) \to  .
 \end{equation}
This implies the first statement because for $n,m$ as stated $H^{(n-1)+(m+1)}_{sing}(W; A) = 0 = H^{n+m+1}_{sing}(W; A)$.

If $m\geq 0$ and $n>d$ or $n+1<-m\leq 0$ then we have $0 = H^{n}_{sing}(W/G) = H^{n,0}(W) = H^{n,1}(W) = \cdots = H^{n,m}(W)$. Similarly if $m\leq 0$ and  $d < n+m\leq n$ or $n< 0$ then we have $H^{n,m}(W) = H^{n,m+1}(W) = \cdots H^{n,-1}(W) = H^{n,0}(W) = H^{n}_{sing}(W/G) = 0$. 
\end{proof}

\begin{theorem}
 If $X$ is a smooth real variety then for any abelian group $A$, there is a strongly convergent spectral sequence 
\begin{equation}\label{ahsstop}
E_{2}^{p,q}(top) \Longrightarrow kr^{p+q}(X;A).
\end{equation}
where $E_{2}^{p,q}(top)= H^{p,-q}(X ; \underline{A})$ for $q\leq 0$ and $E_{2}^{p,q}(top) = 0$ for $q>0$.
\end{theorem}
\begin{proof}
The homotopy fiber sequences (\ref{tophofib}) give the exact couple, 
\begin{equation*}
 D_{2}^{p,q}(top) = \pi_{-p-q,0}(\mcal{K}_{top}^{(-q)}(X) ; A) \;\;\; \text{and} \;\;\; E_{2}^{p,q}(top) = \pi_{-p-q,0}(\mcal{M}_{top}^{(-q)}(X) ; A),
\end{equation*}
which gives rise to a half-plane spectral sequence with exiting differentials in the terminology of \cite[Section 6]{Board:ss}. We have $\pi_{s}\mcal{K}^{(n)}_{top}(\R) = 0 = \pi_{s}\mcal{K}^{(n)}_{top}(\R)^{G}$ for $s< n$ by Theorem \ref{mtow}(5) together with Proposition \ref{srec}. Therefore $\lim_{n}\pi_{s}\mcal{K}^{(n)}_{top}(\R) = 0 = \lim_{n}\pi_{s}\mcal{K}^{(n)}_{top}(\R)^{G}$ for all $s$. These towers satisfy the trivial Mittag-Leffler condition and so $\lim^{1}_{n}\pi_{s}\mcal{K}^{(n)}_{top}(\R) = 0 = \lim^{1}_{n}\pi_{s}\mcal{K}^{(n)}_{top}(\R)^{G}$ for all $s$ (see for example \cite[Section 3.5]{Weibel:hom}). From the Milnor exact sequence \cite[Lemma 5.41]{Thom:ket} we conclude that $\holim \mcal{K}^{(n)}_{top}(\R)$ and $\holim \mcal{K}^{(n)}_{top}(\R)^{G}$ are both weakly contractible. The homotopy limit $\holim \mcal{K}^{(n)}_{top}(\R)$ is computed by first replacing the maps of this tower with $G$-fibrations between fibrant objects and then taking an inverse limit and so we see that $(\holim \mcal{K}^{(n)}_{top}(\R))^{G} \wkeq \holim \mcal{K}^{(n)}_{top}(\R)^{G}$. We conclude that $\holim \mcal{K}^{(n)}_{top}(\R)$ is equivariantly weakly contractible, and therefore $\holim\shom{}{\sing_{\bullet}X(\C)}{\mcal{K}^{(n)}_{top}(\R)}$ is weakly contractible as well. In turn this implies that  $\lim_{p} D_{2}^{p,n-p}(top) = 0$. Finally, because $E_{2}^{p,q}(top) = 0$ when $p > 2\dim X$ by Proposition \ref{bvan}(1) we conclude by  \cite[Theorem 6.1]{Board:ss} that the spectral sequence is strongly convergent and converges to $\colim_{p} D_{2}^{p,n-p}(top) = \pi_{-n,0}(\mcal{K}_{top}^{}(X) ; A)$.
\end{proof}

To summarize, we have three strongly convergent spectral sequences and natural maps between them for a smooth real variety $X$ and abelian group $A$,
\begin{equation}\label{comparess}
 \begin{CD}
  E_{2}^{p,q}(alg)= \H^{p-q}_{\mcal{M}}(X ; A(-q)) \Longrightarrow K_{-p-q}(X;A) \\
@VVV \\
E_{2}^{p,q}(sst) = L^{-q}H\R^{p,-q}(X; A) \Longrightarrow K\R^{sst}_{-p-q}(X;A) \\
@VVV \\
E_{2}^{p,q}(top) = H^{p,-q}(X(\C) ; \underline{A}) \Longrightarrow kr^{p+q}(X(\C);A).
 \end{CD}
\end{equation}

We finish by comparing our construction with Dugger's \cite{dug:kr}. We warn the reader that the indexing conventions used there are different from ours. 

For $n\geq 0$, define $\mcal{K}^{top}_{(n-1)}(\R)$ so that it fits into the homotopy fiber sequence of naive $G$-spectra
\begin{equation}\label{deffib}
\mcal{K}^{(n)}_{top}(\R) \to \mcal{K}_{top}(\R) \to \mcal{K}^{top}_{(n-1)}(\R) .
\end{equation}
There are homotopy fiber sequences of naive $G$-spectra
\begin{equation}\label{hfs}
\mcal{M}^{(n)}_{top}(\R) \to \mcal{K}^{top}_{(n)}(\R) \to \mcal{K}^{top}_{(n-1)}(\R)
\end{equation}
and the towers $\{\mcal{K}^{top}_{(n-1)}(\R)\}$ and $\{\mcal{K}_{top}^{(n)}(\R)\}$ give rise to the same spectral sequence (see for example \cite[Appendix B]{GM:Tate}).

Write $\Omega^{\infty}E$ for the infinite loop space associated to a naive $G$-spectrum $E$. By Proposition \ref{ktopkr} and Proposition \ref{mtopem} we have that $\Omega^{\infty}\mcal{M}^{(n)}_{top}(\R)\wkeq K(\underline{\Z},\R^{n,n})$ and $\Omega^{\infty}\mcal{K}^{top}(\R)\wkeq \Z\times BU$. Taking associated infinite loop spaces of the homotopy fiber sequences (\ref{hfs}) gives the homotopy fiber sequences of $G$-spaces
\begin{equation}\label{dcomphfib}
K(\underline{\mathbb{Z}},\R^{n,n}) \to \Omega^{\infty}\mcal{K}^{top}_{(n)}(\R)\to \Omega^{\infty}\mcal{K}^{top}_{(n-1)}(\R) .
\end{equation}

We briefly recall Dugger's construction, for complete details see \cite[Section 3]{dug:kr}. For a based $G$-space $X$, define $P_{2n}X$ to be the nullification of $X$ at the set of $G$-spaces $\mcal{A}_{\R^{n,n}}\eqdef\{S^{W} \,|\, \R^{n,n}\subsetneq W \}\cup \{S^{W}\wedge G_+  \,|\, \R^{n,n}\subsetneq W \}$ (where the $W$ are $G$-representations). Equivalently, this may be viewed as the Bousfield localization at the set of maps $\{S^{W} \to * \,|\, \R^{n,n}\subsetneq W \}\cup \{S^{W}\wedge G_+ \to * \,|\, \R^{n,n}\subsetneq W \}$. Dugger defines his spectral sequence \cite[Corollary 1.3]{dug:kr} via the tower $\{P_{2n}(\Z\times BU)\}$. 

\begin{proposition}\label{dugagr}
There is a map of towers of $G$-spaces 
$$
\{P_{2n}\Omega^{\infty}\mcal{K}^{top}(\R)\} \to \{\Omega^{\infty}\mcal{K}^{top}_{(n)}(\R)\}
$$
such that $P_{2n}\Omega^{\infty}\mcal{K}^{top}(\R)\to \Omega^{\infty}\mcal{K}^{top}_{(n)}(\R)$ is an equivariant weak equivalence. 
\end{proposition}
\begin{proof}
From the homotopy fiber sequences (\ref{dcomphfib}) and the fact that $\pi_{a,b}K(\underline{\mathbb{Z}},\R^{k,k}) = H^{k-a,k-b}(pt,\underline{\Z})=0$ whenever $a\geq k$,  $b\geq k$ and $a+b > 2k$, we see that if $W$ is a $G$-representation with $\R^{n,n}\subsetneq W$ then  $\pi_{W}\Omega^{\infty}\mcal{K}^{top}_{(n)}(\R) \subseteq  \pi_{W}\Omega^{\infty}\mcal{K}^{top}_{(-1)}(\R) =0$. Similarly if $\R^{n,n}\subsetneq W $, then we also have that $[S^{W}\wedge G_{+},\Omega^{\infty}\mcal{K}^{top}_{(k)}(\R)]_{G} = \pi_{|W|}\Omega^{\infty}\mcal{K}^{top}_{(k)}(\R) = 0$. Therefore by the universal property of the functors $P_{2n}$, the map $\Omega^{\infty}\mcal{K}^{top}(\R) \to \Omega^{\infty}\mcal{K}^{top}_{(n)}(\R)$ factors through the localization $\Omega^{\infty}\mcal{K}^{top}(\R)\to P_{2n}\Omega^{\infty}\mcal{K}^{top}(\R)$. Thus we obtain unique up to homotopy maps $P_{2n}\Omega^{\infty}\mcal{K}^{top}(\R) \to \Omega^{\infty}\mcal{K}^{top}_{(n)}(\R)$ from Dugger's Postnikov tower  so that
\begin{equation*}
\xymatrix@-.2pc{
\Omega^{\infty}\mcal{K}^{top}(\R)\ar[r]\ar[dr] & \Omega^{\infty}\mcal{K}^{top}_{(n)}(\R)\ar[r]& \Omega^{\infty}\mcal{K}_{(n-1)}^{top}(\R)\\
 & P_{2n}\Omega^{\infty}\mcal{K}^{top}(\R) \ar[u]\ar[r] & P_{2n-2}P_{2n}\Omega^{\infty}\mcal{K}^{top}(\R) = P_{2n-2}\Omega^{\infty}\mcal{K}^{top}(\R) \ar[u]
}
\end{equation*}
commutes up to homotopy. Thus we obtain the homotopy commutative comparison diagram of homotopy fiber sequences

\begin{equation*}
\xymatrix{
K(\underline{\mathbb{Z}},\R^{n,n}) \ar[r]\ar[d] & P_{2n}\Omega^{\infty}\mcal{K}^{top}(\R)\ar[r]\ar[d] & P_{2n-2}\Omega^{\infty}\mcal{K}^{top}(\R)\ar[d] \\
K(\underline{\mathbb{Z}},\R^{n,n}) \ar[r] & \Omega^{\infty}\mcal{K}^{top}_{(n)}(\R) \ar[r] & \Omega^{\infty}\mcal{K}^{top}_{(n-1)}(\R) .
}
\end{equation*}
To see that the map $ K(\underline{\mathbb{Z}},\R^{n,n}) \to K(\underline{\mathbb{Z}},\R^{n,n})$ is an equivariant weak equivalence it is enough to verify that it is an isomorphism on $\pi_{n,n}$ and an isomorphism on $\pi_{2n }$ upon forgetting the equivariant structure \cite[Definition 1.4, Lemma 3.7]{Lewis:EM}. It is an isomorphism on $\pi_{2n}$ because the tower $\{\mcal{K}^{(n)}(\Delta_{top}^{\bullet}\times_{\C}\C)\}$ has the nonequivariant homotopy type of the usual Postnikov tower \cite[Theorem 3.4]{FHW:sst}. The forgetful map $\pi_{n,n}K(\underline{\Z}, \R^{n,n}) \to \pi_{2n}K(\underline{\Z}, \R^{n,n})$ is seen to be an isomorphism using the sequence (\ref{forles}) because $H^{0,-1}(pt;\underline{\Z}) = 0 = H^{1,-1}(pt;\underline{\Z})$. We conclude that all the maps $ K(\underline{\mathbb{Z}},\R^{n,n}) \to K(\underline{\mathbb{Z}},\R^{n,n})$ appearing in the comparison diagram above are equivariant weak equivalences.

By induction we see that the induced maps $\pi_{k}P_{2n}\Omega^{\infty}\mcal{K}^{top}(\R)\to \pi_{k}\Omega^{\infty}\mcal{K}^{top}_{(n)}(\R)$ and $\pi_{k}P_{2n}\Omega^{\infty}\mcal{K}^{top}(\R)^{G}\to \pi_{k}\Omega^{\infty}\mcal{K}^{top}_{(n)}(\R)^{G}$ are isomorphisms for all $k\geq 1$. Now we address the case $k=0$. By Theorem \ref{mtow} together with Proposition \ref{srec} we see that $\pi_{s}\mcal{K}^{(n)}_{top}(\R) = 0 = \pi_{s}\mcal{K}^{(n)}_{top}(\R)^{G}$ for $s< n$. This together with the homotopy fiber sequence (\ref{deffib}) of naive $G$-spectra shows that
the maps $\pi_{s}\Omega^{\infty}\mcal{K}^{top}_{(n)}(\R) \to \pi_{s}\Omega^{\infty}\mcal{K}^{top}_{(n-1)}(\R)$ and $\pi_{s}\Omega^{\infty}\mcal{K}^{top}_{(n)}(\R)^{G} \to \pi_{s}\Omega^{\infty}\mcal{K}^{top}_{(n-1)}(\R)^{G}$ are isomorphisms for $s\leq n-1$. Together with the corresponding fact for the tower $P_{2n}\Omega^{\infty}\mcal{K}^{top}(\R)$ (which follows from \cite[Proposition 3.8]{dug:kr}) allows us to conclude the case $k=0$. 
Thus the maps $P_{2n}\Omega^{\infty}\mcal{K}^{top}(\R)\to \Omega^{\infty}\mcal{K}^{top}_{(n)}(\R)$ are equivariant weak equivalences and we are done.
\end{proof}

\begin{remark}
As a consequence the spectral sequence $\{E_{2}^{p,q}(top)\}$ agrees with Dugger's for $p+q \leq 0$, with the caveat that his indexing conventions differ from ours. (Note that Dugger's spectral sequence is fringed and confined to $p+q \leq 0$ while we can have non-zero $E_{2}^{p,q}$ terms for $p+q >0$, as long as $q\leq 0$, because our spectral sequence is constructed using naive $G$-spectra). 
\end{remark}

\section{Conjectures and Computations}\label{conj}
In this section we introduce an equivariant version of Suslin's conjecture on morphic cohomology and we show that the complex, real, and equivariant versions of these conjectures are all equivalent. Based on this we compute real morphic cohomology in weight one, generalizing a partial computation of \cite{Teh:HT}. We also prove a real version of the Bloch-Ogus formula \cite[Corollary 7.4]{BO:gersten}. This allows us to prove that real algebraic equivalence coincides with homological equivalence in codimension two for a large class of real varieties. 
Using these techniques, we compute real and equivariant morphic cohomology of certain real varieties which together with the spectral sequences constructed in the previous section allows us to make new computations of real semi-topological $K$-theory in terms of Atiyah's $KR$-theory. As another application of the spectral sequence we give a simple proof of the mod-$2$ Lichtenbaum-Quillen conjecture over $\R$ generalizing a proof  of Karoubi-Weibel \cite{KW:real}.

Recall Suslin's conjecture on morphic cohomology:
\begin{conjecture} [Suslin's Conjecture in weight $q$]
Let $X$ be a smooth, quasi-projective complex variety. The cycle map
\begin{equation*}
L^{q}H^{p}(X) \to H^{p}_{sing}(X(\C);\Z)
\end{equation*}
is an isomorphism for $p\leq q$ and an injection for $p\leq q+1$.
\end{conjecture}

We introduce the analogues of this conjecture for real morphic cohomology and for equivariant morphic cohomology. Real morphic cohomology is introduced by Friedlander-Walker in \cite{FW:real} (for equivariant morphic cohomology see \cite{HV:VT}).

\begin{conjecture}[Real Suslin's Conjecture in weight $q$]
Let $X$ be a smooth, quasi-projective real variety. The cycle map 
\begin{equation*}
L^{q}H\R^{a+q}(X)=L^{q}H\R^{a,q}(X) \to H^{a,q}(X(\C);\underline{\Z})
\end{equation*}
is an isomorphism for $a \leq 0$ and an injection for $a \leq 1$. 
\end{conjecture} 

\begin{conjecture}[Equivariant Suslin's Conjecture in weight $q$]
Let $X$ be a smooth, quasi-projective real variety. The cycle map 
\begin{equation*}
L^{q}H\R^{a,b}(X) \to H^{a,b}(X(\C);\underline{\Z})
\end{equation*}
is an isomorphism for $a \leq 0$ (and $a\leq q$, $b\leq q$) and an injection for $a=1$ (and $a\leq q$, $b\leq q$). 
\end{conjecture}

\begin{theorem}
\label{equiv}
The following are all equivalent statements:
\begin{enumerate}
\item Suslin's conjecture is valid in weight $q$,
\item Real Suslin's conjecture is valid in weight $q$,
\item Equivariant Suslin's conjecture is valid in weight $q$.
\end{enumerate}
\end{theorem}
\begin{proof}
Recall that if $X$ is a complex variety considered as a real variety, then $\mcal{Z}^{q}(X_{\C}) = \mcal{Z}^{q}(X)\times \mcal{Z}^{q}(X)$, where $G$ acts by switching the factors. Consequently $L^{q}H\R^{a,q}(X) = L^{q}H^{a+q}(X)$ when $X$ is a complex variety viewed as a real variety. Thus if the real Suslin's conjecture is valid in weight $q$ then Suslin's conjecture  is valid as well. The real Suslin's conjecture is a special case of the equivariant Suslin's conjecture. We show that the validity of Suslin's conjecture in weight $q$ implies the real Suslin's conjecture in weight $q$ and that the validity of the real Suslin's conjecture in weight $q$  implies the validity of the equivariant Suslin's conjecture in weight $q$.

Assume that Suslin's conjecture is true in weight $q$. For any $\ell \geq 1$, by \cite[Corollary 5.11]{HV:VT} the map $L^{q}H\R^{a,q}(X;\Z/2^{\ell}) \to H^{a,q}(X(\C);\Z/2^{\ell})$ is an isomorphism for $a\leq 0$ and an injection for $a =1$. Therefore it is enough to verify the real Suslin's conjecture with $\Z[1/2]$-coefficients. For a smooth real variety $X$, let $\pi:X_{\C}\to X$ be the projection and consider the commutative diagram

\begin{equation}
\label{transfer}
 \xymatrix{
 L^{k}H\R^{a,q}(X;\Z[1/2]) \ar[d]\ar[r]^{\pi^{*}} & L^{k}H^{a+q}(X_{\C};\Z[1/2]) \ar[r]^{\pi_{*}}\ar[d] & L^{k}H\R^{a,q}(X;\Z[1/2]) \ar[d] \\
 H^{a,q}(X(\C);\underline{\Z[1/2]}) \ar[r]^{\pi^{*}} & H_{sing}^{a+q}(X(\C);\Z[1/2]) \ar[r]^{\pi_{*}} & H^{a,q}(X(\C);\underline{\Z[1/2]}).
 }
 \end{equation}
The compositions $\pi_{*}\pi^{*} = 2$ and thus the horizontal compositions are isomorphisms. It follows that Suslin's conjecture implies the real Suslin's conjecture.

Now we assume the validity of the real Suslin's conjecture and conclude that the equivariant Suslin's conjecture is valid. Let $F_{q} = \hofib(\mcal{Z}^{q}(X_{\C}) \to \mcal{Z}^{q}_{top}(X(\C)))$ be the homotopy fiber of the cycle map. The homotopy fiber construction is equivariant and so we obtain the equivariant homotopy fiber sequence
\begin{equation*}
F_{q}\to \mcal{Z}^{q}(X_{\C}) \to \mcal{Z}^{q}_{top}(X(\C)).
\end{equation*}
The real Suslin's conjecture implies that $\pi_{k}F_{q}^{G} = 0$ for $k\geq q-1$ and also $\pi_{k}F = 0$ for $k\geq q-1$ (because the real Suslin's conjecture implies Suslin's conjecture). We see therefore that $\Omega^{q-1}F_{q}$ is equivariantly contractible for $q\geq 1$ and for $q=0$, $F_{0}$ is equivariantly contractible as well. It follows that the equivariant Suslin's conjecture holds.  
\end{proof}
\begin{remark}
\label{1/2}
In general, from  Diagram \ref{transfer} we see that if Suslin's conjecture is valid in weight $k$ then 
$$
 L^{k}H\R^{a,b}(X;\Z[1/2]) \to H^{a,b}(X(\C);\underline{\Z[1/2]})
 $$
is an isomorphism for $a+b\leq k$  and an injection for $a+b \leq k+1$. A consequence of this observation is the fact that the cycle map
\begin{equation}\label{mp}
 L^{k}H\R^{a,b}(X;\Z) \to H^{a,b}(X(\C);\underline{\Z})
 \end{equation}
 has 2-torsion kernel  $a+b \leq k+1$ and 2-torsion cokernel for any  $a+b \leq k$, whenever Suslin's conjecture is valid in weight $k$. Moreover if the above isomorphism extends for $a+b\leq k+1$  then the map \ref{mp} has a 2-torsion cokernel for any  $a+b \leq k+1$. This remark will be used below in Theorem \ref{imp}.
\end{remark}

\begin{corollary}
\label{Sw1}
Let $X$ be a smooth, quasi-projective real variety. Then the cycle map
\begin{equation*}
L^{1}H\R^{a,b}(X) \to H^{a,b}(X(\C);\underline{\Z})
\end{equation*}
is an isomorphism for $a\leq 0$ and an monomorphism for $a\leq 1$. 
\end{corollary}
\begin{proof} 
By \cite{FL:algco} Suslin's conjecture holds in weight 1.  According to Theorem \ref{equiv}, this implies that the equivariant Suslin's conjecture holds in weight 1.
\end{proof}

Friedlander-Walker introduce in \cite{FW:real} the equivalence relation of \textit{real algebraic equivalence} for algebraic cycles on a real variety. Briefly two cycles $\alpha$, $\beta$ on a real variety $X$ are real algebraically equivalent provided there is a smooth real curve $C$, two real points $c_{0}$, $c_{1}$ in the same connected component (in the analytic topology) of $C(\R)$, and a cycle $\gamma$ on $X\times C$ such that $\alpha=\gamma|_{c_{0}}$ and $\beta = \gamma|_{c_{1}}$. We write $A^{q}_{rae}(X)$ for the group of codimension-$q$ cycles modulo real algebraic equivalence. If $X$ is smooth then $A^{q} _{rae}(X) = L^{q}H\R^{q,q}(X)$ by \cite[Proposition 4.6]{DS:equiDT} and \cite[Corollary 4.19]{HV:VT}. 

We obtain an improvement of the computation of real Lawson homology of the space of divisors given in \cite{Teh:HT}. 
\begin{corollary} 
\label{div}
Let $X$ be a smooth, quasi-projective real variety. Then
$$
L^{1}H\R^{a,1}(X) \iso H^{a,1}(X(\C);\underline{\Z})\iso H^{a+1} _G(X(\C);\Z(1))
$$
for $a\leq 0$ and 
$$
L^{1}H\R^{1,1}(X)=A^1 _{rae}(X) \hookrightarrow H^{1,1}(X(\C);\underline{\Z})\hookrightarrow H^2 _G(X(\C), \Z(1)).
$$
Here $H^{*}_{G}(X(\C);\Z(q))$ denotes Borel cohomology, and $\Z(q)$ is the sheaf $\Z$ where $G$ acts by $(-1)^{q}$. 

In particular, real algebraic equivalence coincides with homological equivalence for divisors on real smooth varieties.
\end{corollary}
\begin{proof}
 This follows by setting $b=1$ in  Corollary \ref{Sw1}. The relation between Bredon cohomology and Borel cohomology is a particular case of \cite[Proposition 5.17]{HV:VT}. The last statement is just a reformulation of the monomorphism  $L^{1}H\R^{1,1}(X)= A^1 _{rae}(X) \hookrightarrow H^{1,1}(X(\C);\underline{\Z})$ and the equality $L^{1}H\R^{1,1}(X)=A^1 _{rae}(X)$ is proved in \cite{DS:real} and in \cite{HV:VT}.
\end{proof}
\begin{remark} Notice that $L^{1}H\R^{a,1}(X)=0$ for any $a\leq -1$ because $H^r _G(X, \Z(1))=0$ for any $r\leq 0$. In particular, the above computation gives exactly the vanishing range obtained in \cite{Teh:HT}.
 \end{remark}

\begin{conjecture}($2$-adic Lichtenbaum-Quillen conjecture)
Let $X$ be a smooth real variety. For any $k\geq 1$ the map 
$$
\pi_{n}(K(X);\Z/2^{k}) \to \pi_{n}(K(X)^{hG};\Z/2^{k})
$$
is an isomorphism for $n\geq \dim X -1$ and an injection for $n\geq \dim X -2$. 
\end{conjecture}
By \cite[Proposition 4.7, Corollary 3.10]{FW:real} this is equivalent to the statement that for any smooth real variety $X$ and $k\geq 1$, the map
\begin{equation*}
K\R^{sst}_{n}(X;\Z/2^k)\to kr^{-n}(X(\C);\Z/2^{k})
\end{equation*}
is an isomorphism for $n\geq \dim X -1$ and an injection for $n\geq \dim X -2$.

This was conjecture was verified in some cases by Friedlander-Walker \cite{FW:real}, by Karoubi-Weibel \cite{KW:real} (in slightly weaker form for varieties with real points), and in general by Rosenschon-Ostvaer \cite{RO:LQ}. We provide a simple proof of this result using Voevodsky's verification of the Milnor conjecture \cite{Voev:miln} in conjunction with the comparison of spectral sequences (\ref{comparess}). 

\begin{theorem}\label{QL}
 Let $X$ be a smooth, quasi-projective real variety. Then for any $k\geq 1$, 
\begin{equation*}
K\R^{sst}_{n}(X;\Z/2^k)\to kr^{-n}(X(\C);\Z/2^{k})
\end{equation*}
is an isomorphism for $n \geq \dim X -1$ and a monomorphism for $n = \dim X -2$
\end{theorem}
\begin{proof}
We write $E(sst)$ (respectively $E(top)$) for the semi-topological (respectively the topological) spectral sequence constructed in Section \ref{sstss} and \ref{topss}. We show that the map $E_{r}^{p,q}(sst) \to E_{r}^{p,q}(top)$ from (\ref{comparess}) is an isomorphism for $p+q \leq 1 -\dim X $ and a monomorphism for $p+q = 2- \dim X $. Since these spectral sequences are strongly convergent, this will yield the result.

First the case $r=2$. Suppose that $p+q \leq \ 1 - \dim X$ then either $p\leq 0$ or $q \leq -d$. If $p\leq 0$ then by \cite[Corollary 5.10]{HV:VT} 
\begin{equation*}
 E_{2}^{p,q}(sst) = L^{-q}H^{p,-q}(X;A) \to H^{p,-q}(X(\C);\underline{A}) = E_{2}^{p,q}(top)
\end{equation*}
is an isomorphism (note that if both $p,-q \leq 0$ then $E_{2}^{p,q}(sst) = 0 = E_{2}^{p,q}(top)$). If $q\leq -d$ then by Corollary \ref{geqd} then $E_{2}^{p,q}(sst)  \to E_{2}^{p,q}(top)$ is an isomorphism. Similarly if $p+q \leq 2-\dim X$ then either $p\leq 1$ or $q\leq -d$ and we see by \cite[Corollary 5.10]{HV:VT} and Corollary \ref{geqd} that $E_{2}^{p,q}(sst)  \to E_{2}^{p,q}(top)$ is a monomorphism. Proceeding by induction as in \cite[Theorem 6.1]{FHW:sst} shows that $E_{r}^{p,q}(sst)  \to E_{r}^{p,q}(top)$ is an isomorphism for all $r$ and $p+q \leq 1 - \dim X$ and a monomorphism for all $r$ and $p+q\leq 2 - \dim X$.
\end{proof}

\begin{remark}
The above proof is along the same lines as in \cite{KW:real}, a key difference being that there they use a spectral sequence based on the Borel cohomology of $X(\C)$. If $X(\R) \neq \emptyset$ (i.e. the $G$-action on $X(\C)$ is not free) then this spectral sequence will not converge and in this way the use of Bredon cohomology simplifies the situation.  
\end{remark}

If we know that Suslin's conjecture is valid for a particular variety then we can deduce an integral agreement of real semi-topological $K$-theory and $kr$-theory.
\begin{theorem}\label{ZQL}
Suppose that $X$ is a smooth, quasi-projective real variety such that the cycle maps $L^{q}H\R^{a,q}(X) \to H^{a,q}(X(\C);\underline{\Z})$ are isomorphisms for any $q$, any $a\leq 0$ and injections for $a\leq 1$ (in other words $X$ satisfies the real Suslin's conjecture). Then
\begin{equation*}
K\R_{n}^{sst}(X) \to kr^{-n}(X(\C))
\end{equation*}
is an isomorphism for $n\geq \dim X -1$ and a monomorphism for $n\geq \dim X -2$.
\end{theorem}
\begin{proof}
Follows from the comparison of spectral sequences by the same argument as in Theorem \ref{QL}.
\end{proof}

\begin{corollary} 
Let $X$ be a smooth, quasi-projective real surface. Then 
\begin{equation*}
K\R^{sst}_{n}(X;\Z)\to kr^{-n}(X(\C);\Z)
\end{equation*}
is an isomorphism for $n \geq 1$ and a monomorphism for $n = 0$.
\end{corollary}
\begin{proof} This follows from Corollary \ref{div} and Theorem \ref{ZQL}.
\end{proof}

The following corollary was proved in \cite[Corollary 6.9]{FW:real} using different methods. 
 
\begin{corollary}
Let $X$ be a smooth, quasi-projective real curve. Then 
\begin{equation*}
K\R^{sst}_{n}(X;\Z)\to kr^{-n}(X(\C);\Z)
\end{equation*}
is an isomorphism for $n \geq 0$.
\end{corollary}
\begin{proof} It follows again from Corollary \ref{div} and Theorem \ref{ZQL}.
\end{proof}

In certain special examples of real smooth projective varieties, some of the semi-topological invariants can be computed. For example, suppose that $X$ is a smooth, projective geometrical rationally connected threefold defined over real numbers or that $X$ is generic smooth, projective real cubic of dimension $d$ where $3\leq d\leq 8$, $d\neq 7$. Then:
\begin{theorem} 
\label{imp}
Let $X$ be one of the varieties from the above list. Then
\begin{equation*}
L^{q}H\R^{a,b}(X) \to H^{a,b}(X(\C);\underline{\Z})
\end{equation*}
is an isomorphism for $a \leq 1$ (and $b\leq q$). 
\end{theorem} 
\begin{proof} 
By \cite{Voin:1} and \cite{Voin:2}, for  $X$ is as in hypothesis, the complexification $X _\C$ satisfies
$$
L^qH^{a+b}(X _\C)\iso H^{a+b}(X _\C)
$$
for $a+b\leq q +1$.
By \cite[Corollary 5.11]{HV:VT} $L^{q}H\R^{a,q}(X;\Z/2^{\ell}) \to H^{a,q}(X(\C);\Z/2^{\ell})$ is an isomorphism for $a\leq 0$ and an injection for $a =1$ and any $\ell>0$. In particular, we can conclude that the kernel of the cycle map $L^{q}H\R^{a,q}(X;\Z) \to H^{a,q}(X(\C);\underline{\Z})$ is 2-divisible for any $a \leq 1$ and the cokernel is 2-torsion free for any $a \leq 1$ (see \cite[Proposition 3.1]{Voin:1}).  According to the Remark \ref{1/2}, the kernel and cokernel of the cycle map for  $a \leq 1$ and $b\leq q$ are 2-torsion. This implies that the cycle map $L^{q}H\R^{a,q}(X;\Z) \to H^{a,q}(X(\C);\Z)$ is an isomorphism for any $a\leq 1$. 

Recall that for a based $G$-space $W$ and fixed value of $q\geq 0$, there is a long exact sequence relating the non-equivariant and equivariant homotopy groups
\begin{equation}\label{phiex}
 \cdots \to \pi_{p+q+1}W \xrightarrow{\delta^{(q)}} \pi_{p,q+1}W \to \pi_{p,q}W \xrightarrow{\phi^{(q)}} \pi_{p+q}W \xrightarrow{\delta^{(q)}} 
\cdots \to \pi_{0,q}W \xrightarrow{\phi^{(q)}} \pi_{q}W ,
\end{equation}
(the map $\phi^{(q)}$ is the forgetful map). This sequence arises by considering the homotopy cofiber sequence
\begin{equation*}
 \Z/2_+\to S^{0} \to S^{0,1}
\end{equation*}
and the above long-exact sequence is the one induced by the homotopy fiber sequence
\begin{equation*}
 \mapp{S^{0,1}}{W} \to \mapp{S^{0}}{W} \to \mapp{\Z/2_+}{W}.
\end{equation*}
Applying this homotopy fiber sequence to $Z^q(X_{\C})$ and $Z^q _{top}(X(\C))$ we obtain the following diagram
\begin{equation}
\label{cps}
\xymatrix{
L^{q}H\R^{a-1,b} \ar[r]\ar[d] & L^{q}H^{a-1+b} \ar[r]\ar[d] & L^{q}H\R^{a, b-1} \ar[r]\ar[d] & L^{q}H\R^{a,b} \ar[r]\ar[d] & L^{q}H^{a+b} \ar[d] \\
H^{a-1,b} \ar[r] & H^{a-1+b} \ar[r] & H^{a, b-1} \ar[r] & H^{a,b} \ar[r] & H^{a+b}
}
\end{equation}
  The theorem now follows from the above diagram and reverse induction on $b\leq q$.
\end{proof}

\begin{remark}
Theorem \ref{imp} can be extended, using the same ideas,  to real smooth threefolds and fourfolds with zero cycles on the complexification supported on a proper subvariety of dimension less or equal than two.  We also remark that equivariant Suslin's conjecture holds for these varieties and those from Theorem \ref{imp}.
\end{remark}

Our next results will involve the coniveau spectral sequences for real morphic cohomology and for Bredon cohomology. We make use of \cite{CTHK:BO} for these coniveau spectral sequences. Let $Y\subseteq X$ be a closed subvariety, with open complement $U= X- Y$. Let $\mcal{Z}^{q}(X_{\C})_{Y_{\C}} = \hofib(\mcal{Z}^{q}(X_{\C})\to \mcal{Z}^{q}(U_{\C}))$ and define the equivariant morphic cohomology of $X$ with support on $Y$ by $L^{q}H\R^{q-a,q-b}_{Y}(X)= \pi_{a,b}\mcal{Z}^{q}(X_{\C})_{Y_{\C}}$. Similarly, $H^{q-a,q-b}_{Y(\C)}(X(\C);\underline{\Z}) = H^{q-a,q-b}(X(\C), U(\C);\underline{\Z})$ denotes the Bredon cohomology with support on $Y(\C)$. The generalized cycle map extends to a map on cohomology theories with support $\Phi:L^{q}H\R^{s,t}_{Y}(X)\to H^{s,t}_{Y(\C)}(X(\C);\underline{\Z})$. The following lemma makes explicit that both theories satisfy all the properties required in order to make use of the results of \cite{CTHK:BO}. 

\begin{lemma}\label{NisnD} 
Real morphic cohomology and Bredon cohomology are homotopy invariant cohomologies that satisfy Nisnevich excision.
\end{lemma}
\begin{proof}
Both theories are homotopy invariant. That real morphic cohomology satisfies Nisnevich excision follows from the fact that motivic cohomology does. To be explicit, let 
\begin{equation*}
\xymatrix{
V\ar[r]\ar[d] & Y \ar[d] \\
U \ar@{^{(}->}[r] & X
}
\end{equation*} 
be a distinguished Nisnevich square. Let $\mcal{M}^{(q)}(-)$ be the presheaf of spectra on $Sch/\R$ from section \ref{motss} which computes motivic cohomology (so that $\mcal{M}^{(q)}(-\times\Delta^{\bullet}_{top})$ computes real morphic cohomology). The commutative square 
\begin{equation*}
\xymatrix{
\mcal{M}^{(q)}(X\times W)\ar[r]\ar[d] & \mcal{M}^{(q)}(Y\times W) \ar[d] \\
\mcal{M}^{(q)}(U\times W) \ar@{^{(}->}[r] & \mcal{M}^{(q)}(V\times W)
}
\end{equation*} 
is homotopy cartesian whenever $W$ is smooth. Let $\mcal{F}(W)$ and $\mcal{G}(W)$ denote the homotopy fibers of the top and bottom horizontal arrows. By Proposition \ref{srec} $\mcal{F}(\Delta^{\bullet}_{top})\to \mcal{G}(\Delta^{\bullet}_{top})$ is a weak equivalence. By Proposition \ref{hfseq} we see that these are the homotopy fibers of the horizontal arrows in the diagram
 \begin{equation*}
\xymatrix{
\mcal{M}^{(q)}(X\times \Delta^{\bullet}_{top})\ar[r]\ar[d] & \mcal{M}^{(q)}(Y\times \Delta^{\bullet}_{top}) \ar[d] \\
\mcal{M}^{(q)}(U\times \Delta^{\bullet}_{top}) \ar@{^{(}->}[r] & \mcal{M}^{(q)}(V\times \Delta^{\bullet}_{top}) .
}
\end{equation*}
It follows that real morphic cohomology satisfies Nisnevich descent. 

That Bredon cohomology satisfies Nisnevich descent follows from the fact that 
\begin{equation*}
\xymatrix{
V(\C)\ar[r]\ar[d] & Y(\C) \ar[d] \\
U(\C) \ar@{^{(}->}[r] & X(\C)
}
\end{equation*} 
is an equivariant homotopy pushout diagram of $G$-spaces (see for example \cite{DI:hyp}).
\end{proof}

We will use below the following lemma  that proves the local vanishing for equivariant motivic cohomology.

\begin{lemma}\label{avan}
Let $R=\mcal{O}_{X,x_{1},\ldots, x_{n}}$ be a semi-local ring of a smooth real variety $X$ at the points $x_{1},\ldots, x_{n}\in X$. Then $L^{s}H\R^{a,s}(\spec R) = 0$ if $a>0$. 
\end{lemma}
\begin{proof}
Let $\mcal{M}^{(s)}$ be a presheaf computing weight-$s$ motivic cohomology (for example the presheaf from Section \ref{motss}). If $R$ is as in the statement of the lemma, then $\pi_{i}\mcal{M}^{(s)}(\spec R) = H^{2s-i}_{\mcal{M}}(\spec R;\Z(s))=0$ for $i<s$. By Theorem \ref{recog}  we see that the simplicial abelian group $(d\mapsto\pi_{i}\mcal{M}^{(s)}(\spec R\times_{\R} \Delta^{d}_{top}))$ is contractible for all $i<s$. Applying the spectral sequence for a simplicial spectrum \cite[Corollary 4.22]{Jar:et} shows that  $\pi_{i}\mcal{M}\R^{(s)}_{sst}(\spec R)=L^{s}H\R^{s-i,s}(\spec R) = 0$ for $i<s$. 
\end{proof}

We now prove a real analogue of the Bloch-Ogus formula \cite[Corollary 7.4]{BO:gersten}. Namely for a smooth quasi-projective real variety $X$ we show that $H^{p}_{Zar}(X, \mcal{H}^{0,p})\iso A^{p}_{rae}(X)$, where $\mcal{H}^{a,b}$ is the Zariski sheafification of the presheaf $U\mapsto H^{a,b}(U(\C);\underline{\Z})$ and $A^{p}_{rae}(X)$ is the group of codimension-$p$ cycles modulo real algebraic equivalence.

Observe that by \cite[Proposition 5.17]{HV:VT} and \cite[Proposition 1.15]{DSLF:quat}, the sheaf $\mcal{H}^{0,p}$ agrees with the sheafification of the Borel cohomology presheaf $U\mapsto H^{p}_{G}(U(\C);\Z(p))$ (where $\Z(p)$ is the $G$-module $\Z$ with action given by multiplication by $(-1)^{p}$). Thus the following theorem may be seen as an integral extension, for real varieties, of the Bloch-Ogus formula for etale cohomology with finite coefficients \cite[Theorem 7.7]{BO:gersten}, because Borel cohomology with finite coefficients agrees with etale cohomology with finite coefficients for real varieties \cite{Cox:real}.

\begin{theorem}\label{borae}
Let $X$ be a smooth quasi-projective real variety. There is a natural isomorphism
\begin{equation*}
H^{p}_{Zar}(X,\mcal{H}^{0,p}) \iso A^{p}_{rae}(X).
\end{equation*}
\end{theorem}
\begin{proof}
Let $X^{(p)}$ denote the set of points $x\in X$ whose closure has codimension $p$. If $H^*$ is a cohomology theory with supports, define $H^{i}_{x}(X) = \colim_{ U\subseteq X}H^{i}_{\overline{x}\cap U}(U)$ and $H^{i}(k(x)) = \colim_{ U \subseteq \overline{x}}H^{i}(U)$ (where in both colimits, $U$ ranges over nonempty opens). By \cite[Corollary 5.1.11]{CTHK:BO} if $H$ is  homotopy invariant and satisfies Nisnevich excision then for a smooth quasi-projective $X$, the coniveau spectral sequence 
$E_{1}^{p,q} = \oplus_{x\in X^{(p)}}H^{p+q}_{x}(X) \Longrightarrow H^{p+q}(X)$
has $E_{2}$-term given by $H^{p}_{Zar}(X;\mcal{H}^{q})$, where $\mcal{H}^{q}$ is the sheafification of the presheaf $U\mapsto H^{q}(U)$. The associated filtration is $F^{p}H^{n}(X) = \cup_{Z}Im(H^{n}_{Z}(X) \to H(X))$, where the union is over closed subvarieties $Z\subseteq X$ of codimension $p$. 

Equivariant morphic cohomology and Bredon cohomology  satisfy these properties by Lemma \ref{NisnD} and using Lemma \ref{gisin}, the coniveau spectral sequences can be written as
\begin{equation*}
E_{1}^{p,q}(k,b) = \bigoplus\limits_{x\in X^{(p)}}L^{k-p}H\R^{q,b-p}(k(x)) \Longrightarrow L^{k}H\R^{p+q,b}(X)
\end{equation*}
and
\begin{equation*}
\widetilde{E}_{1}^{p,q}(b) = \bigoplus\limits_{x\in X^{(p)}}H^{q,b-p}(k(x);\underline{\Z}) \Longrightarrow H^{p+q,b}(X(\C);\underline{\Z}) ,
\end{equation*}
and the $E_{2}$ terms are respectively given by $E_{2}^{p,q}(k,b) = H_{Zar}^{p}(X;\mcal{L}^{k}\mcal{H}\R^{q,b})$ and $\widetilde{E}_{2}^{p,q}(b)= H_{Zar}^{p}(X;\mcal{H}^{q,b})$.
Observe that $F^{p}L^{p}H\R^{p,p}(X)=L^{p}H\R^{p,p}(X)$ (and $F^{p+1}L^{p}H\R^{p,p}(X)=0$) and so  $E_{\infty}^{p,0}(p,p) = L^{p}H\R^{p,p}(X)$. 

In fact we show that we have $E_{2}^{p,0}(p,p) = L^{p}H\R^{p,p}(X)$. By Lemma \ref{avan} we see that $E_{1}^{p-r,r-1}(p,p) = \oplus_{x\in X^{(p-r)}}L^{r}H\R^{r-1,r}(k(x)) = 0$ for $r>1$. Also $E_{1}^{p+r,-r+1}(p,p) = 0$ if $r\geq 1$. We conclude that $E_{2}^{p,0}(p,p) = E_{\infty}^{p,0}(p,p)$, and therefore $H^{p}_{Zar}(X,\mcal{L}^{p}\mcal{H}\R^{0,p}) = A^{p}_{rae}(X)$.

Consider the commutative diagram 
\begin{equation*}
\xymatrix{
{\bigoplus_{x\in X^{(p-1)}}L^{1}H\R^{0,1}(k(x))} \ar[r]\ar[d] & {\bigoplus_{x\in X^{(p)}}L^{0}H\R^{0,0}(k(x))}\ar[d]\ar[r] & 0 \\
\bigoplus_{x\in X^{(p-1)}}H^{0,1}(k(x);\underline{\Z}) \ar[r] & \bigoplus_{x\in X^{(p)}}H^{0,0}(k(x)) \ar[r] & 0 ,
}
\end{equation*}
the cohomology of the top row computes $H^{p}_{Zar}(X,\mcal{L}^{p}\mcal{H}\R^{0,p})$ and the bottom computes $H^{p}_{Zar}(X,\mcal{H}\R^{0,p})$. The right vertical arrow is an isomorphism and by Corollary \ref{div} the left vertical arrow is also an isomorphism. Therefore $$
A^{p}_{rae}(X)\iso H^{p}_{Zar}(X;\mcal{L}^{p}\mcal{H}\R^{0,p}) \to H^{p}_{Zar}(X;\mcal{H}^{0,p})
$$
 is an isomorphism. 
\end{proof}

We can now prove a real analogue of a result due to Bloch-Srinivas \cite{BS:corr}.
\begin{theorem}\label{cod2}
Let $X$ be a smooth projective real variety of dimension $d$. Suppose that there is a closed subvariety $V\subseteq X _\C$ with $\dim V\leq 2$ such that $CH_{0}(X_{\C}\setminus V) = 0$. Then 
$$
L^{2}H\R^{2,2}(X) = A^{2}_{rae}(X) \hookrightarrow H^{2,2}(X(\C);\underline{\Z}),
$$
 so that real algebraic equivalence on $X$ in codimension 2 agrees with homological equivalence.
\end{theorem}
\begin{proof}
We use the same notations as in the previous theorem for the coniveau spectral sequence for Bredon cohomology. We know that $\widetilde{E}_{\infty}^{2,0}(2) = F^{2}H^{2,2}(X(\C);\underline{\Z})$ because $F^{3}H^{2,2}(X(\C);\underline{\Z}) = 0$.  Because $\widetilde{E}_{1}^{p,q}(2) = 0$ if $q\leq 1$ and $p>2$ or if $p<0$ we see that $\widetilde{E}_{3}^{2,0}(2) = \widetilde{E}_{\infty}^{2,0}(2) = F^{2}H^{2,2}(X(\C);\underline{\Z})\subseteq H^{2,2}(X(\C);\underline{\Z})$. We thus have the exact sequence $\widetilde{E}_{2}^{0,1}(2) \xrightarrow{d_{2}} \widetilde{E}_{2}^{2,0}(2) \to \widetilde{E}_{3}^{2,0}(2)\subseteq H^{2,2}(X(\C);\underline{\Z})$ which is the exact sequence
\begin{equation*}
H^{0}_{Zar}(X, \mcal{H}^{1,2}) \to H^{2}_{Zar}(X,\mcal{H}^{0,2}) \xrightarrow{\alpha} H^{2,2}(X(\C),\underline{\Z}).
\end{equation*}
The cycle map $\Phi$ induces the comparison of coniveau spectral sequences which together with the preceding theorem and its proof gives us the commutative diagram
\begin{equation*}
\xymatrix{
E_{2}^{2,0}(2,2) = H^{2}_{Zar}(X, \mcal{L}^{2}\mcal{H}\R^{0,2}) \ar[r]^{\iso}\ar[d]^{\iso} & E_{\infty}^{2,0}(2,2) = L^{2}H\R^{2,2}(X) \ar[d]^{\Phi} \\
\widetilde{E}_{2}^{2,0}(2) = H^{2}_{Zar}(X, \mcal{H}^{0,2}) \ar@{->>}[r] & \widetilde{E}_{\infty}^{2,0}(2) \subseteq H^{2,2}(X(\C);\underline{\Z}).
}
\end{equation*}
We see that under the identification $H^{2}_{Zar}(X,\mcal{H}^{0,2}) = A^{2}(X)_{rae}$ the map $\alpha$ is identified with the cycle map and so the result will follow from showing that $H^{0}_{Zar}(X, \mcal{H}^{1,2})= 0$. Let $\mcal{H}_{csing}^{k}$ denote the sheaf $U\mapsto H^{k}_{sing}(U(\C);\Z)$. With $\Z[1/2]$-coefficients the composition induced by the map $\pi: X _\C\rightarrow X$ 
$$
H^{0}_{Zar}(X, \mcal{H}^{1,2}[1/2]) \xrightarrow{\pi^{*}} H^{0}_{Zar}(X _\C, (\mcal{H}_{csing}^{3})[1/2]) \xrightarrow{\pi_{*}} H^{0}_{Zar}(X, \mcal{H}^{1,2}[1/2])
$$ 
is an isomorphism. Since $H^{0}_{Zar}(X _\C, (\mcal{H}_{csing}^{3})[1/2])=0$ for a variety $X$ satisfying the condition from hypothesis (by  the proof of \cite[Theorem 1(ii)]{BS:corr}) it is enough to show that $H^{0}_{Zar}(X, \mcal{H}^{1,2})$ is $2$-torsion free. Write $\mcal{L}^{q}\mcal{H}\R^{a,b}_{\Z/2^{k}}$ for the sheafification of $U\mapsto L^{q}H\R^{a,b}(U;\Z/2^{k})$ and similarly for $\mcal{H}^{a,b}_{\Z/2^{k}}$. Consider the comparison of short exact sequences
\begin{equation*}
\xymatrix{
0 \ar[r] & \mcal{L}^{2}\mcal{H}\R^{0,2}\otimes \Z/2^{k} \ar[r]\ar[d] & \mcal{L}^{2}\mcal{H}\R^{0,2}_{\Z/2^{k}} \ar[r]\ar[d] & Tor(\mcal{L}^{2}\mcal{H}\R^{1,2}, \Z/2^{k}) \ar[r]\ar[d] & 0 \\
0 \ar[r] & \mcal{H}^{0,2}\otimes \Z/2^{k} \ar[r] & \mcal{H}\R^{0,2}_{\Z/2^{k}} \ar[r] & Tor(\mcal{H}^{1,2}, \Z/2^{k}) \ar[r] & 0 .
}
\end{equation*}
The map $\mcal{L}^{2}\mcal{H}\R^{0,2}_{\Z/2^{k}} \to \mcal{H}\R^{0,2}_{\Z/2^{k}}$ is an isomorphism by \cite[Corollary 5.13]{HV:VT} and $\mcal{L}^{2}\mcal{H}\R^{1,2} = 0$ by Lemma \ref{avan}. We conclude that $Tor(\mcal{H}^{1,2}, \Z/2^{k}) =0$ and we are done.

\end{proof}
 Because all the real varieties from Theorem \ref{imp} fulfill the condition of Theorem \ref{cod2}, we can use Theorem \ref{cod2} to improve the range of indexes in Theorem \ref{imp}. We will restrict to the case of real geometrically rationally connected varieties of dimension 3, although the same type of arguments and computations can be done for more general varieties that fulfill the condition from Theorem \ref{cod2}.  
 
 \begin{corollary} Let $X$ be a smooth projective real geometrically rationally connected threefold. Then
 \begin{equation*}
L^{q}H\R^{a,b}(X) \to H^{a,b}(X(\C);\underline{\Z})
\end{equation*}
is an isomorphism for $a < q$ (and $b\leq q$) and a monomorphism if $a=q$ (and $b\leq q$). In particular, all equivariant morphic cohomology groups of $X$ are finitely generated.
 \end{corollary}
\begin{proof} 
Write $F  _q = \hofib(\mcal{Z}^{q}(X_{\C}) \to \map{X_{\C}(\C)}{\mcal{Z}_{0}(\A^{q}_{\C})})$ for the homotopy fiber of the cycle map. The homotopy fiber construction is equivariant and yields an equivariant homotopy fiber sequence
\begin{equation*}
 F _q\to\mcal{Z}^{q}(X_{\C}) \to \map{X_{\C}(\C)}{\mcal{Z}_{0}(\A^{q}_{\C})}.
\end{equation*}
Applying Theorem \ref{cod2} and Theorem \ref{imp} we obtain that $\pi _*(F _q^G)=0$. From \cite{Voin:1} we have that $\pi _*(F _q)=0$. We conclude that $F$ is an equivariantly contractible space. Applying the long exact sequences for the above equivariant homotopy fiber sequence we obtain that $\pi _{0,p}(\mcal{Z}^{q}(X_{\C}))\hookrightarrow \pi _{0,p}(\map{X_{\C}(\C)}{\mcal{Z}_{0}(\A^{q}_{\C})})$ for any $p\geq 0$ and that $\pi _{q,p}(\mcal{Z}^{q}(X_{\C}))\simeq \pi _{q,p}(\map{X_{\C}(\C)}{\mcal{Z}_{0}(\A^{q}_{\C})})$ for any $q>0$ and $p\geq 0$. This concludes the statement of the corollary.
\end{proof}
The following corollary computes the real semi-topological K-theory of a real geometrically connected variety.
\begin{corollary} Let $X$ be a smooth projective real geometrically rationally connected threefold. Then
$$
K\R^{sst}_{0}(X; \Z)\hookrightarrow kr^0(X;\Z)
$$
and 
$$
K\R^{sst}_{p}(X; \Z)\xrightarrow{\iso} kr^{-p}(X;\Z)
$$
for any $p\geq 1$.
 \end{corollary}
\begin{proof}
This follows from the previous corollary together with a comparison of spectral sequences constructed in the previous section (see \ref{comparess}). 
\end{proof}

\section{Compatibility of Poincare Dualities}\label{comppd}

The main result of this section asserts that the Poincare duality proved in \cite{HV:VT} between the equivariant morphic cohomology and Dos Santos's equivariant Lawson homology for a smooth, projective real variety is compatible with the Poincare duality in Bredon homology and cohomology of a compact real manifold (see Theorem \ref{poincomp}). The method used here is similar to that of  \cite{FL:dual} where Friedlander-Lawson show that the duality between morphic cohomology and Lawson homology for smooth projective complex varieties is compatible with Poincare duality between singular cohomology and homology. An important consequence of this compatibility is Corollary \ref{geqd} which is proved below and is needed for our applications in section \ref{conj} of the spectral sequence constructed in this paper.

Let $X$ be a smooth real variety of dimension $d$. There are generalized cycle maps $\Phi:L^{q}H\R^{a,b}(X)\to H^{a,b}(X(\C);\underline{\Z})$ relating equivariant morphic cohomology and Bredon cohomology (see \cite{FW:real,HV:VT}). The equivariant Lawson homology defined by dos Santos \cite{DS:real} comes equipped with operations $s:L_{q}H\R_{a,b}(X) \to L_{q-1}H\R_{a,b}(X)$. In particular, iterating the $s$-map together with the equivariant Dold-Thom theorem \cite{DS:equiDT,LF:Gequiv} defines the generalized cycle map $\Psi=s^{q}:L_{q}H\R_{a,b}(X) \to L_{0}H\R_{a,b}(X) = H_{a,b}(X(\C);\underline{\Z})$.

A \textit{real n-bundle} over a $G$-space $W$ is a complex $n$-plane bundle $\pi:E\xrightarrow{} W$  together with an involution $\tau:E\to E$ covering the involution $\sigma:W\to W$. A $G$-manifold is called a \textit{real n-manifold} if  the tangent bundle is a real $n$-bundle over $W$, with involution given by $d\sigma$.
A real $n$-manifold $W$ has a fundamental class $[W]\in H_{n,n}(W;\underline{\Z})$ and satisfies Poincare duality (\cite{DSLF:quat},\cite{LLM}),
\begin{equation*}
 [W]\cap - : H^{p,q}(W;\underline{\Z}) \xrightarrow{\iso} H_{n-p,n-q}(W;\underline{\Z}).
\end{equation*}

The following is the main theorem of this section, whose proof we postpone until the end of the section.
\begin{theorem}\label{poincomp}
Let $X$ be a smooth, projective real variety of dimension $d$. The following square commutes for all $k,l \leq q$, where  the horizontal arrows are the Poincare duality isomorphisms
\begin{equation*} 
\xymatrix{
  L^qH\R^{k,l}(X)  \ar[r]^-{\mcal{D}}\ar[d]_{\Phi} &  L_{d-q}H\R_{d-k,d-l}(X) \ar[d]^{\Psi} \\
H^{k,l}(X(\C),\underline{\Z}) \ar[r]^-{\mcal{P}} & H_{d-k,d-l}(X(\C),\underline{\Z})  .
 }
\end{equation*}
\end{theorem}

We make use of cohomology with supports in the proof of Corollary \ref{geqd} below. The definitions of these are recalled prior to Theorem \ref{borae}.

\begin{corollary}\label{geqd}
Let $U$ be a smooth, quasi-projective real variety of dimension $d$. Then the cycle map
\begin{equation*}
 \Phi:L^{q}H\R^{a,b}(U) \to H^{a,b}(U(\C);\underline{\Z})
\end{equation*}
is an isomorphism for $q\geq d$ (and $a\leq q$, $b\leq q$).
\end{corollary}
\begin{proof}
If $U$ is smooth and projective this follows immediately from Theorem \ref{poincomp} together with the fact that the map $L_{d-q}H\R_{d-a,d-b}(U) \to H_{d-a,d-b}(U(\C),\underline{\Z})$ is an isomorphism for any $q\geq d$. 

Reverse induction on $b$, together with the diagram (\ref{cps}) and the corresponding result for morphic cohomology of complex varieties, shows that in order to prove the general result we may assume that $b=q$. We proceed by induction on $\dim U$ to show that $\Phi$ is an isomorphism whenever $q\geq \dim U$ and $U$ is smooth and quasi-projective. 

The case $\dim U = 0$ is done and so we may assume that the result is valid for all smooth quasi-projective $U$, with $\dim U < d$. Let $U$ be a quasi-projective real variety with $\dim U = d$ and let $U\subseteq X$ be a smooth projective closure with closed complement $Z\subseteq X$ (we may assume that $U$ is connected so that $\dim Z < \dim U$). First we argue that $L^{q}H\R^{a,q}_{W}(X)\to H^{a,q}_{W(\C)}(X(\C);\underline{\Z})$ is an isomorphism for $q\geq d$ and $a\leq q$ whenever $W\subseteq X$ is a proper closed subvariety. Using Lemma \ref{gisin} below we find that it is an isomorphism whenever $W$ is smooth. In particular the case $\dim W = 0$ is done and we proceed by induction on $\dim W$. Let $W_{s}\subseteq W$ be the singular locus. Since $\dim W_{s} < \dim W$ we have by induction that $L^{q}H\R^{a,q}_{W_{s}}(X)\to H^{a,q}_{ W_{s}(\C)}(X(\C))$ is an isomorphism for $q\geq d$ and $a\leq q$. Considering the comparison of long exact sequence 
\begin{equation*}
\xymatrix@-.2pc{
\ar[r] & L^{q}H\R_{W_{s}}^{a,q}(X) \ar[r]
\ar[d] & L^{q}H\R_{W}^{a,q}(X) \ar[r]\ar[d] & L^{q}H\R^{a,q}_{W\setminus W_{s}}(X\setminus W_{s}) \ar[r]\ar[d]  & \\
\ar[r] & H_{W_{s}(\C)}^{a,q}(X(\C)) \ar[r] & H_{W(\C)}^{a,q}(X(\C)) \ar[r] & H^{a,q}_{(W\setminus W_{s})(\C)}((X\setminus W_{s})(\C))  \ar[r] & ,
}
\end{equation*} 
together with Lemma \ref{gisin} below and the inductive hypothesis, we see that for $q\geq d$ and $a\leq q$ the map $L^{q}H\R^{a,q}_{W\setminus W_{s}}(X\setminus W_{s})\to H^{a,q}_{(W\setminus W_{s})(\C)}((X\setminus W_{s})(\C))$ is also an isomorphism. We conclude that
$L^{q}H\R^{a,q}_{W}(X)\to H^{a,q}_{ W(\C)}(X(\C))$ is an isomorphism for $q\geq d$ and $a\leq q$ and any closed, proper subset $W\subseteq X$ (for the case $a=q$ observe that $L^{q}H\R^{a,q}_{W\setminus W_{s}}(X\setminus W_{s}) = 0 = H^{a,q}_{(W\setminus W_{s})(\C)}((X\setminus W_{s})(\C))$ because $q\geq d > \dim W\setminus W_{s}$) .

Now from the comparison of long exact sequences
\begin{equation*}
\xymatrix{
\ar[r] & L^{q}H\R_{Z}^{a,q}(X) \ar[r]
\ar[d] & L^{q}H\R^{a,q}(X) \ar[r]\ar[d] & L^{q}H\R^{a,q}(U) \ar[r]\ar[d]  & \\
\ar[r] & H_{Z(\C)}^{a,q}(X(\C)) \ar[r] & H^{a,q}(X(\C)) \ar[r] & H^{a,q}(U(\C))  \ar[r] & 
}
\end{equation*} 
we see that $L^{q}H\R^{a,q}(U) \to H^{a,q}(U(\C);\underline{\Z})$ is also an isomorphism for $q\geq d$ and $a \leq q$ (for the case $a=q$, note that $L^{q}H\R^{q,q}(X) \to L^{q}H\R^{q,q}(U)$ is surjective and $H_{Z(\C)}^{q+1,q}(X(\C))= 0$ by Proposition \ref{bvan} and the long exact sequence for cohomology with supports, because $2q+1 \geq \dim X(\C)$).
\end{proof}

\begin{lemma}\label{gisin}
Let $X$ be a smooth, quasi-projective real variety and let $Z\subseteq X$ be a smooth closed subvariety of codimension $c$. There  are natural isomorphisms $L^{q-c}H\R^{a-c,b-c}(Z) \xrightarrow{\iso} L^{q}H\R^{a,b}_{Z}(X)$ and $H^{a-c,b-c}(Z(\C);\underline{\Z}) \xrightarrow{\iso} H^{a,b}_{Z(\C)}(X(\C);\underline{\Z})$ which fit into the commutative square
\begin{equation*}
\xymatrix{
L^{q-c}H\R^{a-c,b-c}(Z) \ar[r]^{\iso}\ar[d]^{\Phi} & L^{q}H\R^{a,b}_{Z}(X) \ar[d]^{\Phi} \\
H^{a-c,b-c}(Z(\C);\underline{\Z}) \ar[r]^{\iso}& H^{a,b}_{Z(\C)}(X(\C);\underline{\Z}).
}
\end{equation*}
\end{lemma}
\begin{proof}
Let $N$ denote the normal bundle of $Z\subseteq X$. We consider $Z$ as a closed subvariety of $N$ via the zero-section. The lemma will follow by seeing that we have the commutative diagram
\begin{equation*}
\xymatrix{
L^{q-c}H\R^{a-c,b-c}(Z) \ar[r]^{\iso}\ar[d]^{\Phi} & L^{q}H\R^{a,b}_{Z}(N) \ar[r]^{\iso}\ar[d]^{\Phi} & L^{q}H\R^{a,b}_{Z}(X) \ar[d]^{\Phi} \\
H^{a-c,b-c}(Z(\C);\underline{\Z}) \ar[r]^{\iso} & H^{a,b}_{Z(\C)}(N(\C);\underline{\Z}) \ar[r]^{\iso}& H^{a,b}_{Z(\C)}(X(\C);\underline{\Z}).
}
\end{equation*}
The first horizontal arrows are Thom isomorphisms and the second horizontal arrows arise from the deformation to the normal cone. 

We apply the results of \cite{Panin:oriented}. To put ourselves in the context of this paper we consider the total equivariant morphic cohomology and Bredon cohomology: $\oplus_{q,a,b} L^{q}H\R_{Z}^{a,b}(X)$ and $\oplus_{a,b}H^{a,b}_{Z(\C)}(X(\C);\underline{\Z})$. These are cohomology theories in the sense of \cite[Definition 2.1]{Panin:oriented}, which amounts to saying these are homotopy invariant, satisfy Nisnevich excision, and have localization exact sequences. Both equivariant morphic cohomology and Bredon cohomology clearly satisfy the first and third property, and the excision property is made explicit in Lemma \ref{NisnD}. As shown in \cite[Theorem 2.2]{Panin:oriented} the deformation to the normal cone yields natural isomorphisms $h_{Z}(N)\xrightarrow{\iso} h_{Z}(X)$ for any cohomology theory $h(-)$. The maps involved in this isomorphism preserve the gradings and therefore we obtain the right-hand commutative square in the diagram above.

Now we need to see that $\Phi$ is compatible with Thom isomorphisms (which Panin calls an orientation \cite[Definition 3.1]{Panin:oriented}) in both theories. By \cite[Theorem 3.35, Theorem 3.5]{Panin:oriented} in order for $L^{*}H\R^{*,*}(-)$ to have Thom isomorphisms it is enough to see that it has  a Chern structure \cite[Definition 3.2]{Panin:oriented} which amounts to having first Chern classes for line bundles which satisfy a naturality condition, a projective bundle formula for $X\times \P^{1}$, and vanishing on trivial line bundles. These are easily deduced from the corresponding properties for motivic cohomology (over $\C$ and over $\R$), Theorem \ref{recog}, and the fact that equivariant morphic cohomology may be computed using the complex $z_{equi}(\A^{q}_{\C},0)(X_{\C}\times_{\C}\Delta^{\bullet}_{top})$. Moreover the Thom isomorphism for $E$ is given by cupping with a Thom class $th(E)$ \cite[Lemma 3.33]{Panin:oriented} and keeping track of gradings we see that $th(E)\in L^{n}H\R^{n,n}_{X}(E)$ where $E$ is an $n$-plane bundle.

In \cite[Proposition 1.13]{DSLF:quat} it is verified that Bredon cohomology has Thom isomorphisms for real bundles (in the above sense). In particular Bredon cohomology viewed as a theory on $Sm/\R$ has Thom isomorphisms. To see that the natural transformation $\Phi$ is compatible with the Thom isomorphisms in suffices to see that $\Phi(\mu)\in H^{n,n}_{X(\C)}(E(\C))$ is a Thom class whenever $\mu\in L^{n}H\R^{n,n}_{X}(E)$ is a Thom class. By the description of the Thom class given in \cite[Proposition 1.13]{DSLF:quat}, it suffices to see that for any closed point $i:x\hookrightarrow X$ the restriction $i^{*}\Phi(\mu)$ is a generator of $H^{*,*}_{x}(i^*E(\C);\underline{\Z})$. Note that $i^*E = \A^{n}_{\R}$ or $\A^{n}_{\C}$ depending on whether $k(x) = \R$, or $\C$. In both cases $\Phi:L^{*}H\R_{x}^{*,*}(i^*E) \to H^{*,*}_{x(\C)}(i^{*}E(\C);\underline{\Z})$ is an isomorphism and we are done. 
\end{proof}

Let $X$ be a quasi-projective real variety and $Y$ a projective real variety. The cap-product pairing (\cite[Proposition 1.6]{FW:funcspc}, \cite[Theorem 7.2]{FL:cocyc})
\begin{equation*}
 \mcal{C}_s(X_{\C})\times \Mor{}{X_{\C}}{ \mcal{C}_r(Y_{\C})}  \to \mcal{C}_{r+s}(X_{\C}\times_{\C} Y_{\C})
\end{equation*}
is defined by sending a pair $(f,Z)$ to the graph of the composite function  $Z\to X_{\C} \to \mcal{C}_{r}(Y_{\C})$. This pairing is $G$-equivariant and by \cite[Theorem 2.6]{FW:funcspc} it induces a continuous pairing
\begin{equation*}
\mcal{C}_s(X_{\C})^{an} \times \Mor{}{X_{\C}}{ \mcal{C}_r(Y_{\C})}^{an}  \to \mcal{C}_{r+s}(X_{\C}\times_{\C} Y_{\C})^{an},
\end{equation*}
which in turn induces the continuous equivariant pairing
\begin{equation*}
\cap:  \mcal{Z}_{r}(X_{\C}) \wedge \Mor{}{X_{\C}}{\mcal{Z}_{0}(\A_{\C}^{q})}    \to \mcal{Z}_{r}(X_{\C}\times_{\C} \A_{\C}^{q}),
\end{equation*}
where for typographical convenience here and below we write
\begin{equation*}
\Mor{}{X_{\C}}{\mcal{Z}_{0}(\A_{\C}^{q})}= \frac{\Mor{}{X_{\C}}{\mcal{C}_{0}(\P_{\C}^{q})}^{+}}{\Mor{}{X_{\C}}{\mcal{C}_{0}(\P_{\C}^{q-1})}^{+}}.
\end{equation*}
We also denote by $\cap$ the induced pairing,
\begin{equation*}
\cap:  L_{p}H\R_{s,t}(X) \otimes L^{q}H\R^{k,l}(X)  \to L_{p-q}H\R_{s-k,t-l}(X).
\end{equation*}
When $X$ is irreducible of dimension $d$, the $d$-cycle $X_{\C}\subseteq X_{\C}$ defines the fundamental class
\begin{equation*}
 \eta_{X}\in \pi_{0,0}\mcal{Z}_d(X_{\C})= L_{d}H\R_{d, d}(X_{\C}),
\end{equation*}
and cap-product with the fundamental class gives a map
\begin{equation*}
 \eta_{X}\cap - : L^{q}H\R^{k,l}(X) \to  L_{d-q}H\R_{d-k, d-l}(X),
\end{equation*}
which is the Poincare duality isomorphism. 

For a $G$-space $W$, define the equivariant pairing $\sqcap$ to be the composition

\begin{equation*}
 \sqcap:  \mcal{Z}_{0}(W)\wedge  \mapp{W}{\mcal{Z}_0(S^{d,d})_{0}}  \xrightarrow{}  \mcal{Z}_{0}(W)\wedge \mcal{Z}_{0}(S^{d,d})_{0}\xrightarrow{\mu}  \mcal{Z}_{0}(W\wedge S^{d,d})_{0}
\end{equation*}
where the first arrow is defined by $(x,f) \mapsto (x, f(x))$ and $\mu$ is the map induced by $(x,y)\mapsto x\wedge y$. Together with the suspension isomorphism in Bredon homology this defines a pairing
\begin{equation*}
 \sqcap : H_{s,t}(W;\underline{\Z})\otimes H^{k,l}(W; \underline{\Z}) \to H_{s+ d-k, t+ d-l}(W\wedge S^{d,d};\underline{\Z}) \iso H_{s-k, t-k}(W;\underline{\Z}).
\end{equation*}

\begin{remark}\label{remzero}
The pairing 
$$
\mcal{Z}_{0}(X_{\C})\wedge \Mor{}{X_{\C}}{\mcal{Z}_{0}(\A_{\C}^{q})} \to \mcal{Z}_{0}(X_{\C})\wedge\mapp{X(\C)}{\mcal{Z}_0(S^{q,q})_{0}}
\xrightarrow{\sqcap} \mcal{Z}_{0}(X(\C)_{+}\wedge S^{d,d})_{0}$$
 agrees with the pairings $\cap$ under the isomorphism $\mcal{Z}_{0}(X(\C)_{+}\wedge S^{d,d})_{0} \iso \mcal{Z}_{0}(X_{\C}\times_{\C}\A^{d}_{\C})$.
\end{remark}

The following is merely a rephrasing of \cite[Proposition 1.14]{DSLF:quat}.
\begin{proposition}\cite[Proposition 1.14]{DSLF:quat}
\label{bredcap}
Let $W$ be a smooth compact $G$-manifold of dimension $d$ and let $[W]\in H_{d,d}(W;\underline{\Z})$ denote its fundamental class. The map $[W]\sqcap -$ defined above agrees with cap-product (in Bredon cohomology and homology) with $[W]$ 
$$
 [W]\cap - = [W]\sqcap - :H^{k,l}(W;\underline{\Z}) \to H_{d-k,d-l}(W;\underline{\Z}) .
$$ 
\end{proposition}

We finish our discussion of pairings by showing that the pairings $\cap$ and $\sqcap$ are compatible under the natural transformations $\Phi$ and $\Psi$.
\begin{proposition}\label{capcomp}
 Let $X$ be a smooth projective real variety. The diagram
\begin{equation*}
 \xymatrix@-1pc{
\pi_{0,0}\mcal{Z}_{d}(X_{\C})\otimes \pi_{k,l}\Mor{}{X_{\C}}{\mcal{Z}_{0}(\A_{\C}^{d})}  \ar[d]^-{\Psi\otimes\Phi}\ar[r]^-{\cap}  & \pi_{k,l}\mcal{Z}_{d}(X_{\C}\times_{\C}\A_{\C}^{d}) \ar[d]^-{\Psi}  \\
\pi_{0,0}\Omega^{d,d}\mcal{Z}_{0}(X(\C))\otimes \pi_{k,l}\map{X(\C)}{\mcal{Z}_0(S^{d,d})_{0}}  \ar[r]^-{\sqcap} & \pi_{k,l}\Omega^{d,d}\mcal{Z}_{0}(X(\C)_{+}\wedge S^{d,d})_{0} 
}
\end{equation*}
commutes.
\end{proposition}
\begin{proof}

Consider the following diagram
\begin{equation*}
 \xymatrix{
\P_{\C}^{d}/\P_{\C}^{d-1}\wedge\mcal{Z}_{d}(X_{\C})\wedge\Mor{}{X_{\C}}{\mcal{Z}_0(\A_{\C}^{d})}  \ar[r]^-{1\wedge \;\cap}\ar[d]_{u\wedge 1} &
\P_{\C}^{d}/\P_{\C}^{d-1}\wedge\mcal{Z}_{d}(X_{\C}\times_{\C}\A_{\C}^{d}) \ar[d]^{u\wedge 1} \\
\mcal{Z}_{0}(\A_{\C}^{d})\wedge \mcal{Z}_{d}(X_{\C})\wedge \Mor{}{X_{\C}}{\mcal{Z}_{0}(\A_{\C}^{d})}  \ar[r]^-{1\wedge \;\cap}\ar[d]_{\mu\wedge 1} & \mcal{Z}_{0}(\A_{\C}^{d})\wedge\mcal{Z}_{d}(X_{\C}\times_{\C}\A_{\C}^{d}) \ar[dd]^{\mu}\\
\mcal{Z}_{d}(\A_{\C}^{d}\times_{\C} X_{\C})\wedge \Mor{}{X_{\C}}{\mcal{Z}_{0}(\A_{\C}^{d})} \ar[d]_{1\wedge \pi^*} & \\
\mcal{Z}_{d}(\A_{\C}^{d}\times_{\C} X_{\C})\wedge \Mor{}{\A_{\C}^{d}\times_{\C} X_{\C}}{\mcal{Z}_{0}(\A_{\C}^d)}  \ar[r]^-{\cap} & \mcal{Z}_{d}(\A_{\C}^{d}\times_{\C} X_{\C}\times_{\C}\A_{\C}^{d})  \\
\mcal{Z}_{0}(X_{\C})\wedge \Mor{}{X_{\C}}{\mcal{Z}_{0}(\A_{\C}^{d})}  \ar[u]^-{\pi^*\wedge\pi^*}\ar[d]_{1\wedge\Psi}\ar[r]^-{\cap} & \mcal{Z}_{0}(X_{\C}\times_{\C}\A_{\C}^{d})\ar[d]^{\iso}\ar[u]_-{\pi^*} \\
\mcal{Z}_{0}(X(\C))\wedge \map{X(\C)}{\mcal{Z}_0(S^{d,d})_{0}}  \ar[r]^-{\sqcap} & \mcal{Z}_{0}(X(\C)_{+}\wedge S^{d,d})_{0} .
}
\end{equation*}
Here $u$ is the map induced by sending $x\in \P^{d}_{\C}$ to $x-\infty \in \mcal{Z}_{p}(\P_{\C}^{d})_{0}$ and $\mu$ is induced by $(x,V)\mapsto \{x\}\times V$. The composition of the left vertical maps induces (by taking an adjoint) the map $\Psi\otimes \Phi$ 
and the adjoint of the right vertical compositions induces the map $\Psi = s^{d}$. Thus to establish the proposition it is enough to establish the commutativity of this diagram.

The first square of this  diagram is evidently commutative. Denote the graph of a morphism $g:W_{\C}\to\mcal{C}_{0}(\P_{\C}^k)$ by $\Gamma(g)\in \mcal{Z}_{d}(W_{\C}\times_{\C} \P_{\C}^{k})$. Let $a\in \P_{\C}^{d}$ and $f:X_{\C}\to \mcal{C}_{0}(\P_{\C}^{d})$. The equality of graphs  $\{a\}\times \Gamma(f) = \Gamma(f\pi|_{\{a\}\times X})\in\mcal{Z}_{d}(\A_{\C}^{d}\times_{\C} X_{\C} \times_{\C} \A_{\C}^{d})$ establishes the commutativity of the second square. Similarly for the third square. 
The fourth square is commutative by remark \ref{remzero}.
\end{proof}

The first step towards proving Theorem \ref{poincomp} will be to reduce to the case $q = d$. For this we need to introduce the equivariant analogue of the $s$-map in morphic cohomology.
Consider the commutative diagram
\begin{equation*}
\xymatrix{
 \displaystyle{ \frac{\Mor{}{X_{\C}}{\mcal{C}_{0}(\P_{\C}^{q})}^{+}}{\Mor{}{X_{\C}}{\mcal{C}_{0}(\P_{\C}^{q-1})}^{+}} } \ar[r]^-{\susp}\ar[d]^-{\mcal{D}} & \displaystyle{ \frac{\Mor{}{X_{\C}}{\mcal{C}_{1}(\P_{\C}^{q+1})}^{+}}{\Mor{}{X_{\C}}{\mcal{C}_{1}(\P_{\C}^{q})}^{+}} } \ar[d]^-{\mcal{D}}\\
 \mcal{Z}_{d}(X_{\C}\times_{\C}\A^{q}) \ar[r]^-{\susp} & \mcal{Z}_{d+1}(X_{\C}\times_{\C}\A^{q+1}_{\C})
}
\end{equation*}
where $\susp$ denotes the algebraic suspension map. By \cite[Corollary 4.20]{HV:VT} the vertical maps are equivariant homotopy equivalences, by the real analogue of Lawson's suspension theorem proved by Lam \cite{Lam:t} (see \cite{LLM:real}) the bottom horizontal map is an equivariant homotopy equivalence. We conclude that the the upper horizontal arrow is a homotopy equivalence. Let $\susp^{-1}$ denote a choice of homotopy inverse to this map and consider the following composition


\begin{equation*}
\Mor{}{X_{\C}}{\mcal{Z}_{0}(\A_{\C}^{q})} \wedge\mcal{Z}_{0}(\P^{1}_{\C})_{0} \xrightarrow{\#}   
\Mor{}{X_{\C}}{\mcal{Z}_{1}(\A_{\C}^{q+2})}  \xrightarrow{\susp^{-1}} 
\Mor{}{X_{\C}}{\mcal{Z}_{0}(\A_{\C}^{q+1})}.
\end{equation*}
Define the $s$-map in equivariant morphic cohomology 
$$
s:L^{q}H\R^{a,b}(X) \to L^{q+1}H\R^{a,b}(X)
$$
to be the map induced by this composition together with $\alpha\in\pi_{1,1}\P^{1}_{\C}$ where the map $\alpha:S^{1,1}=\P^{1}_{\C}\to \mcal{Z}_{0}(\P^{1}_{\C})_{0}$ is defined by $x\mapsto \{x\}-\{\infty\}$. Observe that under the duality morphism, the $s$-map in equivariant morphic cohomology is compatible with the $s$-map in dos Santos' equivariant Lawson homology in the sense that we have a commutative diagram
\begin{equation}\label{Dscomp}
\xymatrix{
L^{q}H\R^{a,b}(X) \ar[r]^-{s}\ar[d]^-{\mcal{D}} & L^{q+1}H\R^{a,b}(X) \ar[d]^-{\mcal{D}} \\
L_{d-q}H\R_{d-a,d-b}(X) \ar[r]^-{s} & L_{d-q-q}H\R_{d-a,d-b}(X).
}
\end{equation}

\begin{proposition}\label{scomp}
The transformations $\Phi$ are compatible with the $s$-maps in the sense that
\begin{equation*}
 \xymatrix{
L^{q}H^{p,r}(X) \ar[d]^-{\Phi} \ar[r]^-{s}& L^{q+1}H^{p,r}(X) \ar[ld]^-{\Phi} \\
H^{p,r}(X(\C);\underline{\Z}) & 
}
\end{equation*} 
commutes.
\end{proposition}
\begin{proof}
The space $\mcal{Z}_{k}(\A^{q}_{\C})$ has the equivariant homotopy type of an equivariant Eilenberg-MacLane space $K(\underline{\Z}, \R^{q-k,q-k})$. By \cite{DS:real} the join pairing 
$$
\mcal{Z}_{0}(\A^{p}_{\C})\wedge \mcal{Z}_{0}(\A^{q}_{\C}) \xrightarrow{\#} \mcal{Z}_{1}(\A_{\C}^{p+q+1})\xleftarrow{\susp}\mcal{Z}_{0}(\A_{\C}^{p+q})
$$
 represents the cup product pairing in Bredon cohomology.

Consider the homotopy commutative square
\begin{small}
\begin{equation*}
 \xymatrix@-1pc{
\displaystyle{ \frac{\Mor{}{X_{\C}}{\mcal{C}_{0}(\P_{\C}^{q})}^{+}}{\Mor{}{X_{\C}}{\mcal{C}_{0}(\P_{\C}^{q-1})}^{+}} }\wedge\mcal{Z}_{0}(\P^{1}_{\C})_{0} \ar[d]^-{\Phi}\ar[r]^-{\#} & \displaystyle{ \frac{\Mor{}{X_{\C}}{\mcal{C}_{1}(\P_{\C}^{q+2})}^{+}}{\Mor{}{X_{\C}}{\mcal{C}_{1}(\P_{\C}^{q+1})}^{+}} } \ar[d]^-{\Phi} \ar[r]^-{\susp^{-1}} & \displaystyle{ \frac{\Mor{}{X_{\C}}{\mcal{C}_{0}(\P_{\C}^{q+1})}^{+}}{\Mor{}{X_{\C}}{\mcal{C}_{0}(\P_{\C}^{q})}^{+}} } \ar[d]^-{\Phi} \\
\map{X(\C)}{\mcal{Z}_{0}(\A^{q})}\wedge \mcal{Z}_{0}(\P^{1})_{0} \ar[r]^-{\#} & \map{X(\C)}{\mcal{Z}_{1}(\A^{q+2})} \ar[r]^-{\susp^{-1}} &  \map{X(\C)}{\mcal{Z}_{0}(\A^{q+1})} 
}
\end{equation*}
\end{small}

Let $\alpha:S^{1,1}=\P^{1}_{\C}\to \mcal{Z}_{0}(\P^{1}_{\C})_{0}$ be the map defined by $x\mapsto \{x\}-\{\infty\}$. The map $\alpha$ represents the canonical generator of $\pi_{1,1}\mcal{Z}_{0}(\P^{1})_{0}\iso \Z$. The upper horizontal composition restricted to $\alpha$ induces the $s$-map on equivariant morphic cohomology. On the other hand the bottom horizontal composition restricted to $\alpha$ induces the identity because it represents the product with  $\alpha\in H^{0,0}(pt;\underline{\Z})= \Z$ and $\alpha$ is the canonical generator of this group.
\end{proof}

The Poincare dualities $\mcal{D}$ and $\mcal{P}$ are both given by capping with a fundamental class. We verify that the transformation $\Psi$ takes fundamental classes in equivariant Lawson homology to fundamental classes in Bredon homology. The cycle map $cyc:\mcal{Z}_d(X) \to H_{d,d}(X(\C);\underline{\Z})$ sends an irreducible closed subvariety $i:V\subseteq X$ to $cyc(V) = i_*[V(\C)]\in H_{d,d}(X(\C);\underline{\Z})$. 

\begin{lemma}\label{cycmap}
 Let $X$ be a projective real variety. The composition of $s$-maps together with the Dold-Thom isomorphism 
\begin{equation*}
\Psi:L_{k}H\R_{k,k}(X) \to H_{k,k}(X(\C);\underline{\Z}), 
\end{equation*}
coincides with the cycle-map $V\mapsto [V(\C)]$.
\end{lemma}
\begin{proof}
 The proof is the similar to \cite[Proposition 6.4]{FM:filt}. Let $V\subseteq X$ be a closed 
irreducible subvariety of dimension $k$ and $i:\tilde{V}\to V\subseteq X$ be a resolution of singularities. Then $i_{*}\eta_{\tilde{V}} = [V_{\C}\subseteq X_{\C}]\in L_{k}H\R_{k,k}(X)$. By Noether normalization 
there is a dominant map $f: V\to \P^{k}_{\R}$ of some degree $d$. Consider the commutative 
diagram
\begin{equation*}
 \xymatrix{
L_kH\R_{k,k}(\P^{k}) \ar[d]_{\iso}^{\Psi} & L_kH\R_{k,k}(\tilde{V}) \ar[l]_{f_*}\ar[r]^{i_*} 
\ar[d]^{\Psi} & L_{k}H\R_{k,k}(X) \ar[d]^{\Psi} \\
H_{k,k}(\P^{k};\underline{\Z}) & H_{k,k}(\tilde{V}(\C);\underline{\Z}) 
\ar[l]_{f_*}\ar[r]^{i_*} & H_{k,k}(X(\C);\underline{\Z}) .
}
\end{equation*}
The map $\Psi: L_kH\R_{k,k}(\P^{k})\to H_{k,k}(\P^{k}(\C);\underline{\Z})$ is an isomorphism \cite{DS:real}, $\Psi(\eta_{\P^{k}}) = [\P^{k}(\C)]$,  and $f_*(\eta_{V}) = d\cdot\eta_{\P^{k}}\in L_kH\R_{k,k}(\P^{k})$. Note that $H^{0,0}(\tilde{V}(\C);\underline{\Z}) \iso \Z$ because 
$$
H^{0,0}(\tilde{V}(\C);\underline{\Z}) \iso 
H^{0}_{sing}(\tilde{V}(\C)/G;\Z)\iso\Hom{}{H_{0}(\tilde{V}(\C)/G;\Z)}{\Z}\iso\Z
$$
and so by equivariant Poincare duality  we have $H_{k,k}(\tilde{V}(\C);\underline{\Z}) \iso \Z$ generated by the fundamental class $[\tilde{V}(\C)]$.

If we see that $f_*:H_{k,k}(\tilde{V}(\C);\underline{\Z}) \to H_{k,k}(\P^k(\C);\underline{\Z})$ is 
multiplication by $d$ we are done. Indeed if $f_*$ is multiplication by $d$ then it follows that 
$\Psi$  takes the natural generator $\eta_{\tilde{V}}$ of $L_kH_{k,k}(\tilde{V}(\C);\underline{\Z})\iso \Z$ to 
the generator $[\tilde{V}(\C)]$ of $H_{k,k}(\tilde{V}(\C);\underline{\Z})\iso \Z$.

The fact that this map is multiplication by $d$ follows from the commutative square
\begin{equation*}
 \xymatrix{
H_{k,k}(\tilde{V}(\C);\underline{\Z}) \ar[r]\ar[d] & H_{2k}(\tilde{V}(\C);\underline{\Z}) 
\ar[d]^{\cdot d}\\
H_{k,k}(\P^{k}(\C);\underline{\Z}) \ar[r]^{\iso} & H_{2k}(\P^{k}(\C);\underline{\Z}),
}
\end{equation*}
where the horizontal maps are the forgetful maps and send the equivariant fundamental class to 
the non-equivariant fundamental class.
\end{proof}

We also need to know that the homotopy invariance in equivariant Lawson homology and the suspension isomorphism in Bredon homology are suitably related. First recall that the suspension isomorphism is given by the following explicit map of spaces.
\begin{lemma}(\cite[Remark 2.14]{LF:Gequiv})
The suspension isomorphism $\phi$ is induced by the adjoint of the map
\begin{equation*}
 S^{d,d}\wedge\mcal{Z}_{0}(W) \xrightarrow{u\times 1} \mcal{Z}_{0}(S^{d,d})_{0}\wedge \mcal{Z}_{0}(W) \xrightarrow{\mu} \mcal{Z}_{0}(S^{d,d}\wedge W_{+})_{0} \xrightarrow{\tau} \mcal{Z}_{0}(W_{+}\wedge S^{d,d})_{0},
\end{equation*}
where the first map is the  map $u(x) = \{x\}-\{\infty\}$ .
\end{lemma}
 
\begin{proposition}\label{suspcomm}
Let $X$ be a smooth, projective real variety. Then 
\begin{equation*}
 \xymatrix{
\mcal{Z}_d(X_{\C}\times_{\C}\A_{\C}^d) \ar[d]^-{\Psi} & \mcal{Z}_{0}(X_{\C}) \ar[l]_-{\pi^{*}} \ar[ld]^-{\phi} \\
\Omega^{d,d}\mcal{Z}_{0}(X(\C)_{+}\wedge S^{d,d})_{0} & 
}
\end{equation*}
is homotopy commutative.
 \end{proposition}
\begin{proof}
By adjointness, the homotopy commutativity of the triangle follows from the homotopy commutativity of the following diagram:
\begin{equation*}
 \xymatrix{
\P_{\C}^{d}/\P_{\C}^{d-1}\wedge\mcal{Z}_{d}(X_{\C}\times_{\C}\A_{\C}^{d}) \ar[d]^{u\times 1} & \P_{\C}^{d}/\P_{\C}^{d-1}\wedge \mcal{Z}_{0}(X_{\C}) \ar[l]_-{1\wedge\pi^*}\ar[d]^-{u\times 1} \\
\mcal{Z}_{0}(\A_{\C}^{d})\wedge\mcal{Z}_{d}(X_{\C}\times_{\C} \A_{\C}^{d}) \ar[d]^-{\mu} & \mcal{Z}_{0}(\A_{\C}^{d})\wedge \mcal{Z}_{0}(X_{\C}) \ar[d]^-{\mu}\ar[l]_-{1\wedge\pi^*} \\
\mcal{Z}_{d}(\A_{\C}^{d}\times_{\C} X_{\C} \times \A_{\C}^{d})  & \mcal{Z}_{0}(\A_{\C}^{d}\times_{\C} X_{\C}) \ar[l]_-{\pi_1^*}\ar[dl]^{\tau^*}\\
\mcal{Z}_{0}(X_{\C}\times_{\C} \A_{\C}^{d})\ar[u]^-{\pi_2^*}  & 
}
\end{equation*}
Here $\pi_1:\A_{\C}^{d}\times_{\C} X_{\C} \times_{\C} \A_{\C}^{d} \to \A_{\C}^{d}\times_{\C} X_{\C}$ and $\pi_2:\A_{\C}^{d}\times_{\C} X_{\C} \times_{\C} \A_{\C}^{d} \to X_{\C}\times_{\C} \A_{\C}^{d}$ are the projections, $\mu$ is multiplication, and $\tau$ is the twist isomorphism. The left vertical composition induces the $s$-map and  the right vertical composition induces the suspension isomorphism $\phi$. The first two squares evidently commute. The two maps 
$ \tau\pi_1$ and $\pi_2:\A_{\C}^{d}\times_{\C} X_{\C} \times_{\C} \A_{\C}^{d} \to \A_{\C}^d\times_{\C} X_{\C}$ 
are homotopic via 
$$F:\A_{\C}^{d}\times_{\C} X_{\C} \times_{\C} \A_{\C}^{d}\times_{\C} \A_{\C}^1 \to \A_{\C}^{d}\times_{\C} X_{\C} ,$$ 
where $F(a,x,b,t) = (ta + (1-t)b,x)$,
and therefore $\pi_2^*\tau^*$ and $\pi_1^*$  are homotopic. 
\end{proof}
We now prove the compatibility 
\begin{varthm}[\textbf{Theorem \ref{poincomp}}]
Let $X$ be a smooth, projective real variety of dimension $d$. The following square commutes, where  the horizontal arrows are the Poincare duality isomorphisms
\begin{equation*} 
\xymatrix{
  L^qH\R^{k,l}(X)  \ar[r]^-{\mcal{D}}\ar[d]_{\Phi} &  L_{d-q}H\R_{d-k,d-l}(X) \ar[d]^{\Psi} \\
H^{k,l}(X(\C),\underline{\Z}) \ar[r]^-{\mcal{P}} & H_{d-k,d-l}(X(\C),\underline{\Z})  .
 }
\end{equation*}
\end{varthm}
\begin{proof}
Proposition together with the commutativity of diagram (\ref{Dscomp}) shows that in order to prove the theorem it suffices to consider the case $q=d$ and to prove that the following square commutes
\begin{equation*} 
\xymatrix{
  L^dH\R^{d-k,d-l}(X) = \pi_{k,l}\mcal{Z}^d(X_{\C}) \ar[r]^-{\mcal{D}}\ar[d]_{\Phi} & \pi_{k,l}\mcal{Z}_{0}(X_{\C}) = L_{0}H\R_{k,l}(X) \ar[d]^{\Psi} \\
H^{d-k,d-l}(X(\C),\underline{\Z}) \ar[r]^-{\mcal{P}} & H_{k,l}(X(\C),\underline{\Z}) .
 }
\end{equation*}
In this case the map $\Psi$ is simply the Dold-Thom isomorphism.

Consider the following diagram 

\begin{equation*}\label{pdiag}
 \xymatrix{
 \pi_{k,l}\Mor{}{X_{\C}}{\mcal{Z}_{0}(\A_{\C}^{d})}  \ar[d]^-{\Phi}\ar[r]^-{\eta_{X}\cap -}  & \pi_{k,l}\mcal{Z}_{d}(X_{\C}\times_{\C}\A_{\C}^{d}) \ar[d]^-{\Psi} & \pi_{k,l}\mcal{Z}_{0}(X_{\C}) \ar[l]_-{\iso}\ar[d]^-{=} \\
 \pi_{k,l}\map{X(\C)}{\mcal{Z}_0(S^{d,d})_{0}}  \ar[r]^-{[X(\C)]\sqcap -} & \pi_{k,l}\Omega^{d,d}\mcal{Z}_{0}(X(\C)_{+}\wedge S^{d,d})_{0} & \pi_{k,l}\mcal{Z}_{0}(X(\C)) \ar[l]_-{\phi} .
}
\end{equation*}

 The upper-composition is the 
Poincare duality isomorphism between equivariant morphic cohomology and equivariant Lawson homology and by Proposition \ref{bredcap} the composition along the bottom
 row  is the Poincare duality isomorphism between Bredon cohomology and homology.

The first square commutes by Proposition \ref{capcomp} together with Lemma \ref{cycmap} and the second square commutes by Proposition \ref{suspcomm}
\end{proof}

\appendix
\section{$\Gamma$-spaces}\label{gammaspc}
We recall Segal's notion of $\Gamma$-space, introduced by Segal in \cite{Segal:cat}, where he showed that such objects give rise infinite loop spaces. Let $\underline{n}$ be the pointed set $\underline{n}=\{0,1,2,\ldots, n\}$, pointed by $0$. The category $\Gamma^{op}$ is the category whose objects are $\underline{n}$ for $n\geq 0$ with pointed set maps ($\Gamma^{op}$ is equivalent to the opposite of the category $\Gamma$ considered by Segal). A $\Gamma$-object in a category $\mcal{C}$ is a functor $F:\Gamma^{op}\to \mcal{C}$. Generally $\mcal{C}$ will be a pointed category in which case we additionally require that $F(0) = *$. When $\mcal{C}$ is the category of spaces or of simplicial sets then we refer to a $\Gamma$-object $F$ as a \textit{$\Gamma$-space}. A $\Gamma$-space is said said to be \textit{special} if 
\begin{equation*}
 (p_{1*}, \ldots, p_{n*}):F(\underline{n}) \to F(\underline{1})^{\times n}
\end{equation*}
is a weak equivalence, where $p_{i}:\underline{n}\to \underline{1}$ is the function $p_{i}(i) = 1$ and $p_{i}(j) = 0$ if $i\neq j$. If $\mcal{F}$ is a presheaf of special $\Gamma$-spaces then we will say that $\mcal{F}$ is special.

Given a special $\Gamma$-space $F$ one can associate a spectrum, which we denote by $\mathbb{S}F$. First note that a $\Gamma$-space $F$ is in particular a bisimplicial set (or simplicial space) $n\mapsto F(\underline{n})$.  Write $BF$ for the classifying $\Gamma$-space defined by
$$
BF(\underline{m})= |n\mapsto F(\underline{m}\wedge\underline{n}) | ,
$$
here $\underline{m}\wedge \underline{n}$ is identified with $\underline{mn}$ via $(i,j) \mapsto n(i-1) +j$. For a $\Gamma$-space $A$ there is a natural map $S^{1}\wedge A(\underline{1}) \to BA(\underline{1})$. The spectrum associated to $F$ is obtained by evaluating iterated classifying spaces of $F$:
\begin{equation*}
\mathbb{S}F = (F(\underline{1}),\, BF(\underline{1}),\, B^{2}F(\underline{1}),\, \ldots ).
\end{equation*}
By \cite{Segal:cat} the fact that $F$ is special implies that $BF$ is special, that $B^{k}F(\underline{1}) \to \Omega B^{k}F(\underline{1})$ is a weak equivalence for $k\geq 1$, and that $F(\underline{1}) \to \Omega BF(\underline{1})$ is a homotopy group completion.

Now suppose that $G= \Z/2$ acts on the $\Gamma$-space $F$ and that both $\underline{m}\mapsto F(\underline{m})$ and $\underline{m}\mapsto F(\underline{m})^{G}$ are special. Say that such a $\Gamma$-space is \textit{equivariantly special}. Realization of a simplicial space commutes with taking fixed points and therefore we immediately conclude the following proposition.

\begin{proposition}\label{egam}
Let $F$ by an equivariantly special $\Gamma$-space and $\mathbb{S}F$ the associated $\Omega$-spectrum. Then the naive fixed point spectrum $(\mathbb{S}F)^{G}$ is the $\Omega$-spectrum associated to the fixed point $\Gamma$-space $\underline{n}\mapsto F(\underline{n})^{G}$.
\end{proposition}

\section{Globalization}\label{glob}
We recall the Godement resolution which gives an explicit model for computing the Zariski hypercohomology of $X$ with coefficients in a presheaf of spectra \cite{Thom:ket}. A particularly useful property of this resolution is that it behaves appropriately with respect to the Galois action of $Gal(\C/\R)$. 

Suppose that $F(-)$ is a presheaf of fibrant spectra on $Sch/\R$. The Godement complex $F(-) \to G^{\bullet}F(-)$ is defined as follows. Let
\begin{equation*}
GF(U) = \prod_{u\in U}F(\spec \mcal{O}_{U,u})
\end{equation*}
then define $G^{n}F =G\circ\cdots \circ GF$ to be the $(n+1)$-fold composition of $G$. The inclusions and projections of the various factors gives rise to the cosimplicial spectrum $n\mapsto G^{n}F(U)$ and define $\mcal{G}F(U) = \tot G^{\bullet}F(U)$. 

In general for a presheaf of spectra, $F(-)$, define 
\begin{equation*}
\mcal{G}F(U) = \tot G^{\bullet}QF(U)
\end{equation*} 
where $Q$ is a fixed functorial fibrant replacement  for spectra. Then $QF(-)$ is a presheaf of fibrant  spectra and the homotopy groups of $\mcal{G}F(U)$ define the Zariski hypercohomology with coefficients in $F$,
\begin{equation*}
\pi_{k}\mcal{G}F(U) = \mathbb{H}^{-k}_{Zar}(U;F).
\end{equation*}
Moreover $\pi_{k}\mcal{G}F(U) = [S^{k}\wedge U_{+}, F]_{Zar}$, where $[-,-]_{Zar}$ denotes maps in the Zariski-local homotopy category of presheaves of  spectra on $Sch/\R$ (\cite{Jar:spre,Mitch:thom}).

Now we consider the case when $F'$ is a presheaf of  naive $\Z/2$-spectra on $Sch/\R$ obtained by $F' = p_{*}p^{*}F$ for some presheaf $F$ on $Sch/\R$, where $p:Sch/\R \to Sch/\C$ is the map of sites specified by $X\mapsto X_{\C}$).
Observe that if $F$ is additive, i.e $F(U\coprod V) = F(U)\times F(V)$, then for a real variety $U$ we have
\begin{equation*}
GF'(U) = \prod_{u\in U}F((\spec\mcal{O}_{U,u})_{\C}) = \prod_{v\in U_{\C}}F((\spec\mcal{O}_{U_{\C},v})= p_*p^{*}GF(U).
\end{equation*}
Indeed $(\spec \mcal{O}_{U,u})_{\C} = \spec \mcal{O}_{U_{\C},\pi^{-1}\{u\}} = \coprod_{u'_{i}\in\pi^{-1}\{u\}}\spec \mcal{O}_{U_{\C},u'_{i}}$, where $\pi:U_{\C}\to U$ is the projection. 
Take the fibrant replacement functor $Q$ on spectra used above to be a $\Z/2$-equivariant fibrant replacement functor. Further suppose that canonical map $F(U) \to F(U_{\C})^{\Z/2}$ is a weak equivalence. For a cosimplicial spectrum $A^{\bullet}$ we have that $\tot(A^{\bullet})_{n} = \tot(A^{\bullet}_{n})$ and from the description of $\tot(A^{\bullet}_{n})$ as the simplicial function complex $\tot(A^{\bullet}_{n}) = \shom{}{\Delta}{A^{\bullet}_{n}}$ on cosimplicial simplicial sets, we see that $\tot(-)$ commutes with fixed points. Because $(GF(U_{\C}))^{\Z/2} = G(F(U_{\C})^{\Z/2})= GF(U)$ and $[\tot G^{\bullet}F(U)]^{\Z/2} = \tot(G^{\bullet}F(U)^{\Z/2})$ we see that
\begin{equation*}
\mcal{G}F(U)\to(\mcal{G}F(U_{\C}))^{\Z/2} 
\end{equation*} 
is a weak equivalence as well.

The reason for using this specific fibrant replacement functor is the following proposition which is needed during our discussion of the topological spectral sequence in Section \ref{topss}. We remind the reader that the presheaf $\mcal{K}^{(n)}(-)$ is obtained by applying the Godement resolution to the presheaf $\mcal{W}^{(n)}(-) = \mcal{K}_{geom}(-,\P^{\wedge n})$ and $\mcal{K}_{geom}(X,Y)$ is defined to be the spectrum associated to a certain special $\Gamma$-space (Definition \ref{defkgeom}). The presheaf $\mcal{M}^{(n)}$ is obtained by applying the Godement resolution to the presheaf of spectra associated to the presheaf of simplical abelian groups $z_{equi}(\P^{\wedge n}, 0)(-\times\Delta^{\bullet})$. See  Section \ref{motss} for more details on these constructions and Section \ref{sstss} for a recollection of the meaning of the notation $F(\Delta^{\bullet}_{top}\times_{\R}-)$ and $F(\Delta^{\bullet}_{top}\times_{\C}-)$ for a presheaf $F(-)$ on $Sch/\R$, used below. 
\begin{proposition}\label{gres}
Let $X$ be a real variety. Then 
$$
\mcal{K}^{(n)}(\Delta^{\bullet}_{top}\times_{\R} X)\to \mcal{K}^{(n)}(\Delta^{\bullet}_{top}\times_{\C}X)^{\Z/2} 
$$ 
and 
$$
\mcal{M}^{(n)}(\Delta^{\bullet}_{top}\times_{\R} X) \to \mcal{M}^{(n)}(\Delta^{\bullet}_{top}\times_{\C}X)^{\Z/2} 
$$ 
are weak equivalences.
\end{proposition}
\begin{proof}
The spectrum $\mcal{K}_{geom}(X,Y)$ is the spectrum associated to the special $\Gamma$-space $m\mapsto \Hom{\R}{X}{G_{Y}^{(m)}}$. By Proposition \ref{egam} we see that $\mcal{K}_{geom}(X_{\C},Y)^{G} = \mcal{K}_{geom}(X,Y)$ because $\Hom{\R}{X}{G_{Y}^{(m)}}=\Hom{\R}{X_{\C}}{G_{Y}^{(m)}}^{G}$. In Section \ref{motss} the spectrum $\mcal{W}^{(n)}(X)$ is defined by a  homotopy colimit of a certain diagram of spectra, $\mcal{W}^{(n)}(X) = \hocolim_{I} \mcal{K}_{geom}(X, \P_{I}^{\wedge n})$. Taking fixed points commutes with homotopy colimits of naive $G$-spectra (see \cite[Remark 5.6]{DI:hyp} for an argument in terms of $G$-spaces, note that the same argument applies to $G$-simplicial sets and homotopy colimits of spectra are computed degreewise) and so we have that $\mcal{W}^{(n)}(X)\wkeq \mcal{W}^{(n)}(X_{\C})^{G} $. Similarly we see that $\mcal{M}(z_{equi}(\P^{\wedge n}, 0)(X\times_{\R} \Delta_{\R}^{\bullet}))\wkeq \mcal{M}(z_{equi}(\P^{\wedge n}, 0)(X_{C}\times_{\R} \Delta_{\R}^{\bullet}))^{G}$.
These are additive presheaves and so from the discussion above on Godement resolutions we conclude that $\mcal{K}^{(n)}(-)=\mcal{G}\mcal{W}^{(n)}(-)$ and $\mcal{M}^{(n)}(-)= \mcal{G}\mcal{M}(z_{equi}(\P^{\wedge n}, 0)(-\times_{\R} \Delta_{\R}^{\bullet}))$ also satisfy $\mcal{K}^{(n)}(X)\wkeq\mcal{K}^{(n)}(X_{\C})^{\Z/2}$ and $\mcal{M}^{(n)}(X)\wkeq\mcal{M}^{(n)}(X_{\C})^{\Z/2} $. 

For a presheaf $F$, the colimit $F(\Delta^{n}_{top}\times_{\C} X) = \colim_{T\to U(\R)}F(U_{\C}\times_{\R} X)$ is a filtered colimit and so fixed points commutes with this colimit, 
\begin{equation*}
F(\Delta^{n}_{top}\times_{\C} X)^{\Z/2} = \colim_{T\to U(\R)}F(U_{\C}\times_{\R} X)^{\Z/2}.
\end{equation*}
We conclude that $
\mcal{K}^{(n)}(\Delta^{\bullet}_{top}\times_{\R} X)\wkeq \mcal{K}^{(n)}(\Delta^{\bullet}_{top}\times_{\C}X)^{\Z/2}$ and 
$\mcal{M}^{(n)}(\Delta^{\bullet}_{top}\times_{\R} X) \wkeq \mcal{M}^{(n)}(\Delta^{\bullet}_{top}\times_{\C}X)^{\Z/2}$.
\end{proof}

\bibliographystyle{amsalpha} 
\bibliography{spseq}

\providecommand{\bysame}{\leavevmode\hbox to3em{\hrulefill}\thinspace}
\providecommand{\MR}{\relax\ifhmode\unskip\space\fi MR }
\providecommand{\MRhref}[2]{%
  \href{http://www.ams.org/mathscinet-getitem?mr=#1}{#2}
}
\providecommand{\href}[2]{#2}
\begin{thebibliography}{LMSM86}

\bibitem[Ati66]{A:kreal}
M.~F. Atiyah, \emph{{$K$}-theory and reality}, Quart. J. Math. Oxford Ser. (2)
  \textbf{17} (1966), 367--386. \MR{MR0206940 (34 \#6756)}

\bibitem[BG73]{BG}
Kenneth~S. Brown and Stephen~M. Gersten, \emph{Algebraic {$K$}-theory as
  generalized sheaf cohomology}, Algebraic {K}-theory, {I}: {H}igher
  {K}-theories ({P}roc. {C}onf., {B}attelle {M}emorial {I}nst., {S}eattle,
  {W}ash., 1972), Springer, Berlin, 1973, pp.~266--292. Lecture Notes in Math.,
  Vol. 341. \MR{MR0347943 (50 \#442)}

\bibitem[BO74]{BO:gersten}
Spencer Bloch and Arthur Ogus, \emph{Gersten's conjecture and the homology of
  schemes}, Ann. Sci. \'Ecole Norm. Sup. (4) \textbf{7} (1974), 181--201
  (1975). \MR{MR0412191 (54 \#318)}

\bibitem[Boa99]{Board:ss}
J.~Michael Boardman, \emph{Conditionally convergent spectral sequences},
  Homotopy invariant algebraic structures ({B}altimore, {MD}, 1998), Contemp.
  Math., vol. 239, Amer. Math. Soc., Providence, RI, 1999, pp.~49--84.
  \MR{1718076 (2000m:55024)}

\bibitem[BS83]{BS:corr}
S.~Bloch and V.~Srinivas, \emph{Remarks on correspondences and algebraic
  cycles}, Amer. J. Math. \textbf{105} (1983), no.~5, 1235--1253. \MR{714776
  (85i:14002)}

\bibitem[Cox79]{Cox:real}
David~A. Cox, \emph{The \'etale homotopy type of varieties over {${\bf R}$}},
  Proc. Amer. Math. Soc. \textbf{76} (1979), no.~1, 17--22. \MR{MR534381
  (80f:14009)}

\bibitem[CTHK97]{CTHK:BO}
Jean-Louis Colliot-Th{\'e}l{\`e}ne, Raymond~T. Hoobler, and Bruno Kahn,
  \emph{The {B}loch-{O}gus-{G}abber theorem}, Algebraic $K$-theory (Toronto,
  ON, 1996), Fields Inst. Commun., vol.~16, Amer. Math. Soc., Providence, RI,
  1997, pp.~31--94. \MR{MR1466971 (98j:14021)}

\bibitem[DI04]{DI:hyp}
Daniel Dugger and Daniel~C. Isaksen, \emph{Topological hypercovers and
  {$\mathbb A\sp 1$}-realizations}, Math. Z. \textbf{246} (2004), no.~4,
  667--689. \MR{MR2045835 (2005d:55026)}

\bibitem[dS03a]{DS:real}
Pedro~F. dos Santos, \emph{Algebraic cycles on real varieties and {${\mathbb
  Z}/2$}-equivariant homotopy theory}, Proc. London Math. Soc. (3) \textbf{86}
  (2003), no.~2, 513--544. \MR{MR1971161 (2004c:55026)}

\bibitem[dS03b]{DS:equiDT}
\bysame, \emph{A note on the equivariant {D}old-{T}hom theorem}, J. Pure Appl.
  Algebra \textbf{183} (2003), no.~1-3, 299--312. \MR{MR1992051 (2004b:55021)}

\bibitem[dSLF04]{DSLF:quat}
Pedro~F. dos Santos and Paulo Lima-Filho, \emph{Quaternionic algebraic cycles
  and reality}, Trans. Amer. Math. Soc. \textbf{356} (2004), no.~12, 4701--4736
  (electronic). \MR{MR2084395 (2005m:55018)}

\bibitem[Dug05]{dug:kr}
Daniel Dugger, \emph{An {A}tiyah-{H}irzebruch spectral sequence for
  {$KR$}-theory}, $K$-Theory \textbf{35} (2005), no.~3-4, 213--256 (2006).
  \MR{MR2240234 (2007g:19004)}

\bibitem[FHW04]{FHW:sst}
Eric~M. Friedlander, Christian Haesemeyer, and Mark~E. Walker,
  \emph{Techniques, computations, and conjectures for semi-topological
  {$K$}-theory}, Math. Ann. \textbf{330} (2004), no.~4, 759--807.
  \MR{MR2102312}

\bibitem[FL92a]{FL:algco}
Eric~M. Friedlander and H.~Blaine Lawson, Jr., \emph{A theory of algebraic
  cocycles}, Ann. of Math. (2) \textbf{136} (1992), no.~2, 361--428.
  \MR{MR1185123 (93g:14013)}

\bibitem[FL92b]{FL:cocyc}
\bysame, \emph{A theory of algebraic cocycles}, Ann. of Math. (2) \textbf{136}
  (1992), no.~2, 361--428. \MR{MR1185123 (93g:14013)}

\bibitem[FL97]{FL:dual}
Eric~M. Friedlander and H.~Blaine Lawson, \emph{Duality relating spaces of
  algebraic cocycles and cycles}, Topology \textbf{36} (1997), no.~2, 533--565.
  \MR{MR1415605 (97k:14007)}

\bibitem[FM94]{FM:filt}
Eric~M. Friedlander and Barry Mazur, \emph{Filtrations on the homology of
  algebraic varieties}, Mem. Amer. Math. Soc. \textbf{110} (1994), no.~529,
  x+110, With an appendix by Daniel Quillen. \MR{MR1211371 (95a:14023)}

\bibitem[FS02]{FS:AHSS}
Eric~M. Friedlander and Andrei Suslin, \emph{The spectral sequence relating
  algebraic {$K$}-theory to motivic cohomology}, Ann. Sci. \'Ecole Norm. Sup.
  (4) \textbf{35} (2002), no.~6, 773--875. \MR{MR1949356 (2004b:19006)}

\bibitem[FV00]{FV:biv}
Eric~M. Friedlander and Vladimir Voevodsky, \emph{Bivariant cycle cohomology},
  Cycles, transfers, and motivic homology theories, Ann. of Math. Stud., vol.
  143, Princeton Univ. Press, Princeton, NJ, 2000, pp.~138--187. \MR{MR1764201}

\bibitem[FW01a]{FW:compK}
Eric~M. Friedlander and Mark~E. Walker, \emph{Comparing {$K$}-theories for
  complex varieties}, Amer. J. Math. \textbf{123} (2001), no.~5, 779--810.
  \MR{MR1854111 (2002i:19004)}

\bibitem[FW01b]{FW:funcspc}
\bysame, \emph{Function spaces and continuous algebraic pairings for
  varieties}, Compositio Math. \textbf{125} (2001), no.~1, 69--110.
  \MR{MR1818058 (2001m:14033)}

\bibitem[FW02a]{FW:real}
\bysame, \emph{Semi-topological {$K$}-theory of real varieties}, Algebra,
  arithmetic and geometry, Part I, II (Mumbai, 2000), Tata Inst. Fund. Res.
  Stud. Math., vol.~16, Tata Inst. Fund. Res., Bombay, 2002, pp.~219--326.
  \MR{MR1940670 (2003h:19005)}

\bibitem[FW02b]{FW:sstfct}
\bysame, \emph{Semi-topological {$K$}-theory using function complexes},
  Topology \textbf{41} (2002), no.~3, 591--644. \MR{MR1910042 (2003g:19005)}

\bibitem[FW03]{FW:ratisos}
\bysame, \emph{Rational isomorphisms between {$K$}-theories and cohomology
  theories}, Invent. Math. \textbf{154} (2003), no.~1, 1--61. \MR{MR2004456
  (2004j:19002)}

\bibitem[GM95]{GM:Tate}
J.~P.~C. Greenlees and J.~P. May, \emph{Generalized {T}ate cohomology}, Mem.
  Amer. Math. Soc. \textbf{113} (1995), no.~543, viii+178. \MR{MR1230773
  (96e:55006)}

\bibitem[Gra95]{Gray:wt}
Daniel~R. Grayson, \emph{Weight filtrations via commuting automorphisms},
  $K$-Theory \textbf{9} (1995), no.~2, 139--172. \MR{MR1340843 (96h:19001)}

\bibitem[GW00]{GW:geom}
Daniel~R. Grayson and Mark~E. Walker, \emph{Geometric models for algebraic
  {$K$}-theory}, $K$-Theory \textbf{20} (2000), no.~4, 311--330, Special issues
  dedicated to Daniel Quillen on the occasion of his sixtieth birthday, Part
  IV. \MR{MR1803641 (2001m:19006)}

\bibitem[HV09]{HV:VT}
Jeremiah Heller and Mircea Voineagu, \emph{Vanishing theorems for real
  algebraic cycles}, Submitted. Preprint available at
  http://arxiv.org/abs/0909.0569 (2009).

\bibitem[Jar87]{Jar:spre}
J.~F. Jardine, \emph{Simplicial presheaves}, J. Pure Appl. Algebra \textbf{47}
  (1987), no.~1, 35--87. \MR{MR906403 (88j:18005)}

\bibitem[Jar97]{Jar:et}
\bysame, \emph{Generalized \'etale cohomology theories}, Progress in
  Mathematics, vol. 146, Birkh\"auser Verlag, Basel, 1997. \MR{MR1437604
  (98c:55013)}

\bibitem[KW03]{KW:real}
Max Karoubi and Charles Weibel, \emph{Algebraic and {R}eal {$K$}-theory of real
  varieties}, Topology \textbf{42} (2003), no.~4, 715--742. \MR{MR1958527
  (2004c:19004)}

\bibitem[Lam90]{Lam:t}
T.-K. Lam, \emph{Spaces of real algebraic cycles and homotopy theory}, Ph.D.
  thesis, SUNY Stony Brook, 1990.

\bibitem[Lev08]{Levine:ss}
Marc Levine, \emph{The homotopy coniveau tower}, J. Topol. \textbf{1} (2008),
  no.~1, 217--267. \MR{2365658 (2008j:14013)}

\bibitem[Lew92]{Lewis:EM}
L.~Gaunce Lewis, Jr., \emph{Equivariant {E}ilenberg-{M}ac {L}ane spaces and the
  equivariant {S}eifert-van {K}ampen and suspension theorems}, Topology Appl.
  \textbf{48} (1992), no.~1, 25--61. \MR{1195124 (93i:55016)}

\bibitem[LF97]{LF:Gequiv}
P.~Lima-Filho, \emph{On the equivariant homotopy of free abelian groups on
  {$G$}-spaces and {$G$}-spectra}, Math. Z. \textbf{224} (1997), no.~4,
  567--601. \MR{MR1452050 (98i:55014)}

\bibitem[LLFM03]{LLM:real}
H.~Blaine Lawson, Paulo Lima-Filho, and Marie-Louise Michelsohn,
  \emph{Algebraic cycles and the classical groups. {I}. {R}eal cycles},
  Topology \textbf{42} (2003), no.~2, 467--506. \MR{MR1941445 (2003m:14013)}

\bibitem[LMSM86]{LLM}
L.~G. Lewis, Jr., J.~P. May, M.~Steinberger, and J.~E. McClure,
  \emph{Equivariant stable homotopy theory}, Lecture Notes in Mathematics, vol.
  1213, Springer-Verlag, Berlin, 1986, With contributions by J. E. McClure.
  \MR{MR866482 (88e:55002)}

\bibitem[May96]{May:equi}
J.~P. May, \emph{Equivariant homotopy and cohomology theory}, CBMS Regional
  Conference Series in Mathematics, vol.~91, Published for the Conference Board
  of the Mathematical Sciences, Washington, DC, 1996, With contributions by M.
  Cole, G. Comeza{\~n}a, S. Costenoble, A. D. Elmendorf, J. P. C. Greenlees, L.
  G. Lewis, Jr., R. J. Piacenza, G. Triantafillou, and S. Waner. \MR{MR1413302
  (97k:55016)}

\bibitem[Mit97]{Mitch:thom}
Stephen~A. Mitchell, \emph{Hypercohomology spectra and {T}homason's descent
  theorem}, Algebraic {$K$}-theory ({T}oronto, {ON}, 1996), Fields Inst.
  Commun., vol.~16, Amer. Math. Soc., Providence, RI, 1997, pp.~221--277.
  \MR{MR1466977 (99f:19002)}

\bibitem[MVW06]{MVW:mot}
Carlo Mazza, Vladimir Voevodsky, and Charles Weibel, \emph{Lecture notes on
  motivic cohomology}, Clay Mathematics Monographs, vol.~2, American
  Mathematical Society, Providence, RI, 2006. \MR{MR2242284 (2007e:14035)}

\bibitem[Pan03]{Panin:oriented}
I.~Panin, \emph{Oriented cohomology theories of algebraic varieties},
  $K$-Theory \textbf{30} (2003), no.~3, 265--314, Special issue in honor of
  Hyman Bass on his seventieth birthday. Part III. \MR{2064242 (2005f:14043)}

\bibitem[R{\O}05]{RO:LQ}
Andreas Rosenschon and Paul~Arne {\O}stv{\ae}r, \emph{The homotopy limit
  problem for two-primary algebraic {$K$}-theory}, Topology \textbf{44} (2005),
  no.~6, 1159--1179. \MR{MR2168573 (2006d:19001)}

\bibitem[Seg74]{Segal:cat}
Graeme Segal, \emph{Categories and cohomology theories}, Topology \textbf{13}
  (1974), 293--312. \MR{MR0353298 (50 \#5782)}

\bibitem[Sus03]{Suslin:GraySS}
A.~Suslin, \emph{On the {G}rayson spectral sequence}, Tr. Mat. Inst. Steklova
  \textbf{241} (2003), no.~Teor. Chisel, Algebra i Algebr. Geom., 218--253.
  \MR{MR2024054 (2005g:14043)}

\bibitem[SV00]{SV:rel}
Andrei Suslin and Vladimir Voevodsky, \emph{Relative cycles and {C}how
  sheaves}, Cycles, transfers, and motivic homology theories, Ann. of Math.
  Stud., vol. 143, Princeton Univ. Press, Princeton, NJ, 2000, pp.~10--86.
  \MR{MR1764199}

\bibitem[Teh08]{Teh:HT}
Jyh-Haur Teh, \emph{Harnack-{T}hom theorem for higher cycle groups and {P}icard
  varieties}, Trans. Amer. Math. Soc. \textbf{360} (2008), no.~6, 3263--3285.
  \MR{MR2379796}

\bibitem[Tho85]{Thom:ket}
R.~W. Thomason, \emph{Algebraic {$K$}-theory and \'etale cohomology}, Ann. Sci.
  \'Ecole Norm. Sup. (4) \textbf{18} (1985), no.~3, 437--552. \MR{MR826102
  (87k:14016)}

\bibitem[Voe03]{Voev:miln}
Vladimir Voevodsky, \emph{Motivic cohomology with {${\bf Z}/2$}-coefficients},
  Publ. Math. Inst. Hautes \'Etudes Sci. (2003), no.~98, 59--104. \MR{MR2031199
  (2005b:14038b)}

\bibitem[Voi08]{Voin:1}
Mircea Voineagu, \emph{Semi-topological {$K$}-theory for certain projective
  varieties}, J. Pure Appl. Algebra \textbf{212} (2008), no.~8, 1960--1983.
  \MR{MR2414696 (2009b:19003)}

\bibitem[Voi09]{Voin:2}
\bysame, \emph{{Cylindrical Homomorphisms and Lawson Homology}}, To appear in
  Journal of $K$-theory (2009).

\bibitem[Wal96]{Walk:thes}
Mark Walker, \emph{Motivic complexes and the $k$-theory of automorphisms},
  Ph.D. thesis, UIUC, 1996.

\bibitem[Wal00]{Walker:adams}
Mark~E. Walker, \emph{Adams operations for bivariant {$K$}-theory and a
  filtration using projective lines}, $K$-Theory \textbf{21} (2000), no.~2,
  101--140. \MR{MR1804538 (2002i:19005)}

\bibitem[Wal02]{Walk:Thom}
\bysame, \emph{Semi-topological {$K$}-homology and {T}homason's theorem},
  $K$-Theory \textbf{26} (2002), no.~3, 207--286. \MR{MR1936355 (2003k:19004)}

\bibitem[Wei94]{Weibel:hom}
Charles~A. Weibel, \emph{An introduction to homological algebra}, Cambridge
  Studies in Advanced Mathematics, vol.~38, Cambridge University Press,
  Cambridge, 1994. \MR{1269324 (95f:18001)}

\end{thebibliography}
\end{document}